\documentclass{article}
\usepackage[latin9]{inputenc}
\usepackage{mathtools}
\usepackage{amsmath}
\usepackage{amsthm}
\usepackage{amssymb}

\makeatletter
  \theoremstyle{definition}
  \newtheorem{defn}{\protect\definitionname}
  \theoremstyle{remark}
  \newtheorem{rem}{\protect\remarkname}
\theoremstyle{plain}
\newtheorem{thm}{\protect\theoremname}
  \theoremstyle{plain}
  \newtheorem{cor}{\protect\corollaryname}
  \theoremstyle{remark}
  \newtheorem*{rem*}{\protect\remarkname}
  \theoremstyle{plain}
  \newtheorem{prop}{\protect\propositionname}
  \theoremstyle{plain}
  \newtheorem*{thm*}{\protect\theoremname}

\usepackage{amsfonts}
\usepackage{graphicx}
\providecommand{\U}[1]{\protect\rule{.1in}{.1in}}

  \providecommand{\definitionname}{Definition}
  \providecommand{\propositionname}{Proposition}
  \providecommand{\remarkname}{Remark}
  \providecommand{\theoremname}{Theorem}
\providecommand{\corollaryname}{Corollary}
\providecommand{\theoremname}{Theorem}

  \providecommand{\definitionname}{Definition}
  
  \providecommand{\propositionname}{Proposition}
  \providecommand{\remarkname}{Remark}
  \providecommand{\theoremname}{Theorem}
\providecommand{\corollaryname}{Corollary}
\providecommand{\theoremname}{Theorem}

  \providecommand{\definitionname}{Definition}
  \providecommand{\propositionname}{Proposition}
  \providecommand{\remarkname}{Remark}
  \providecommand{\theoremname}{Theorem}
\providecommand{\corollaryname}{Corollary}
\providecommand{\theoremname}{Theorem}

  \providecommand{\definitionname}{Definition}
  \providecommand{\propositionname}{Proposition}
  \providecommand{\remarkname}{Remark}
  \providecommand{\theoremname}{Theorem}
\providecommand{\corollaryname}{Corollary}
\providecommand{\theoremname}{Theorem}

\makeatother

  \providecommand{\definitionname}{Definition}
  \providecommand{\propositionname}{Proposition}
  \providecommand{\remarkname}{Remark}
  \providecommand{\theoremname}{Theorem}
\providecommand{\corollaryname}{Corollary}
\providecommand{\theoremname}{Theorem}

\begin{document}

\title{Extended Hamilton-Jacobi theory, contact manifolds \\
 and integrability by quadratures}

\author{Sergio \ Grillo\\
 {\small{}{}{}{}{}Instituto Balseiro, Universidad Nacional de
Cuyo and CONICET}\\[-5pt] {\small{}{}{}{}{}Av. Bustillo 9500,
San Carlos de Bariloche}\\[-5pt] {\small{}{}{}{}{}R8402AGP, República
Argentina}\\[-5pt] {\small{}{}{}{}{}sergiog@cab.cnea.gov.ar}\\[5pt]
Edith \ Padrón\\[-5pt]{\small{}{}{}{}{}ULL-CSIC Geometría Diferencial
y Mecánica Geométrica}\\[-5pt] {\small{}{}{}{}{}Departamento
de Matemáticas, Estad{í}stica e IO}\\[-5pt] {\small{}{}{}{}{}Universidad
de la Laguna, La Laguna, Tenerife, Canary Islands, Spain}\\[-5pt]
{\small{}{}{}{}{}mepadron@ull.edu.es}}
\maketitle
\begin{abstract}
A Hamilton-Jacobi theory for general dynamical systems, defined on
fibered phase spaces, has been recently developed. In this paper we
shall apply such a theory to contact Hamiltonian systems, as those
appearing in thermodynamics and on geodesic flows in fluid mechanics.
We first study the partial and complete solutions of the Hamilton-Jacobi
Equation (HJE) related to these systems. Then we show that, for a
given contact system, the knowledge of what we have called a complete
\textit{pseudo-isotropic solution} ensures the integrability by quadratures
of its equations of motion. This extends to contact manifolds a recent
result obtained in the context of general symplectic and Poisson manifolds. 
\end{abstract}

\section{Introduction}

In the last few years, different versions of the Hamilton-Jacobi Theory
have been developed beyond the scenario of Hamiltonian systems in
a cotangent bundle. See for instance \cite{bmmp,pepin-holo,pepin-noholo,lmm,hjp}.
In a recent paper \cite{gp}, an extension of the theory to general
dynamical systems (on fibered phase spaces), which contains as particular
cases the previously mentioned versions, has been presented. Let us
briefly describe it. Consider a finite-dimensional smooth manifold
$M$ and a vector field $X\in\mathfrak{X}\left(M\right)$. Suppose
that $M$ is a fibered bundle with base manifold $N$ and surjective
submersion $\Pi:M\rightarrow N$ (\emph{ipso facto} an open map).
Related to these data (see \cite{gp}) we have the $\Pi$-\textbf{Hamilton-Jacobi
equation }($\Pi$-\textbf{HJE}) for $X$: 
\begin{equation}
\sigma_{\ast}\circ\Pi_{\ast}\circ X\circ\sigma=X\circ\sigma,\label{hjer}
\end{equation}
whose unknown is a section $\sigma:N\rightarrow M$ of $\Pi$ (\emph{ipso
facto} a closed map). If $\sigma$ solves the equation above we shall
say that $\sigma$ is a \textbf{(global) solution of the }$\Pi$-\textbf{HJE
for} $X$. On the other hand, given an open subset $U\subseteq M$,
we shall call \textbf{local solution of} \textbf{the }$\Pi$-\textbf{HJE
for} $X$ \textbf{along} $U$ to any solution of the $\Pi_{\left|U\right.}$-HJE
for $X_{\left|U\right.}$. (Here, we are seeing $\Pi_{\left|U\right.}$
as a fibration onto $\Pi\left(U\right)$ and $X_{\left|U\right.}$
as a vector field on $U$). Note that $\sigma$ is a solution of the
$\Pi$-HJE for $X$ if and only if 
\begin{equation}
\sigma_{*}\circ X^{\sigma}=X\circ\sigma,\label{hjrel}
\end{equation}
where $X^{\sigma}\coloneqq\Pi_{*}\circ X\circ\sigma$, i.e. the vector
fields $X\in\mathfrak{X}\left(M\right)$ and $X^{\sigma}\in\mathfrak{X}\left(N\right)$
are $\sigma$-related. (Moreover, it can be shown that $\sigma$ is
a solution of above equation if and only if its image is an $X$-invariant
submanifold). This means that, in order to find the trajectories of
$X$ along the image of $\sigma$, we can look for the integral curves
of $X^{\sigma}$.

The standard Hamilton-Jacobi Theory is obtained when $M$ is a co-tangent
bundle $T^{*}Q$, the fibration $\Pi$ is the canonical projection
$\pi_{Q}:T^{*}Q\rightarrow Q$ and $\sigma:Q\rightarrow T^{*}Q$ is
an exact $1$-form. A Hamilton-Jacobi Theory for closed $1$-forms
may be found, for instance, in Ref. \cite{pepin-holo}.

Given another manifold $\Lambda$ such that $\dim\Lambda+\dim N=\dim M$,
a \textbf{complete solution of the }$\Pi$-\textbf{HJE for }$X$ is
a surjective local diffeomorphism $\Sigma:N\times\Lambda\rightarrow M$
such that, for all $\lambda\in\Lambda$, 
\begin{equation}
\sigma_{\lambda}:=\Sigma\left(\cdot,\lambda\right):p\in N\longmapsto\Sigma\left(p,\lambda\right)\in M\label{psr}
\end{equation}
is a solution of the $\Pi$-HJE for $X$. The local version is obtained
by replacing $M$, $X$, $\Pi$ and $N$ by $U$, $X_{\left|U\right.}$,
$\Pi_{\left|U\right.}$ and $\Pi\left(U\right)$, respectively, being
$U$ an open subset of $M$. Each section $\sigma_{\lambda}$ is called
a \textbf{partial solution}. We showed in \cite{gp} that a (local)
complete solution $\Sigma$ exists around every point $p\in M$ for
which $X\left(p\right)\notin\ker\Pi_{\ast,p}$. Denoting by $\mathfrak{p}_{N}$
and $\mathfrak{p}_{\Lambda}$ the projections of $N\times\Lambda$
onto $N$ and $\Lambda$, respectively, it is easy to prove that a
surjective local diffeomorphism $\Sigma$ is a complete solution if
and only if 
\begin{equation}
\Pi\circ\Sigma=\mathfrak{p}_{N}\ \ \ \text{and}\ \ \ \Sigma_{\ast}\circ X^{\Sigma}=X\circ\Sigma,\label{Srelr}
\end{equation}
being $X^{\Sigma}\in\mathfrak{X}\left(N\times\Lambda\right)$ the
unique vector field on $N\times\Lambda$ satisfying 
\begin{equation}
\left(\mathfrak{p}_{N}\right)_{\ast}\circ X^{\Sigma}=\Pi_{\ast}\circ X\circ\Sigma\ \ \ \text{and}\ \ \ \left(\mathfrak{p}_{\Lambda}\right)_{\ast}\circ X^{\Sigma}=0.\label{Srelrr}
\end{equation}
Note that $X^{\Sigma}\left(p,\lambda\right)=\left(X^{\sigma_{\lambda}}\left(p\right),0\right)$,
with $X^{\sigma_{\lambda}}:=\Pi_{\ast}\circ X\circ\sigma_{\lambda}\in\mathfrak{X}\left(N\right)$,
so, in particular, 
\begin{equation}
\mathsf{Im}X^{\Sigma}\subseteq TN\times\left\{ 0\right\} .\label{xsin}
\end{equation}
Also, the fields $X$ and $X^{\Sigma}$ are $\Sigma$-related. This
implies that all the trajectories of $X$ can be obtained from those
of $X^{\Sigma}$. More precisely, since each trajectory of $X^{\Sigma}$
is clearly of the form $t\mapsto\left(\gamma\left(t\right),\lambda\right)\in N\times\Lambda$,
for some fixed point $\lambda\in\Lambda$ {[}see Eq. \eqref{xsin}{]},
those of $X$ are given by $t\mapsto\Sigma\left(\gamma\left(t\right),\lambda\right)$.
So, for each $\lambda$, it rests to find the curves $\gamma$, which
are the integral curves of the vector field $X^{\sigma_{\lambda}}\in\mathfrak{X}\left(N\right)$.
In the standard case, it is well-known that such curves can be found
up to quadratures. This gives rise to a natural problem: in the general
case, when can such curves be constructed up to quadratures? Or, which
conditions must a complete solution satisfy, in the general case,
in order to ensure integrability by quadratures of the vector fields
$X^{\sigma_{\lambda}}$?

When $M$ is a general symplectic or a Poisson manifold and $X$ is
a Hamiltonian vector field, we have shown in \cite{gp} that, under
a rather simple additional condition, the integral curves of each
$X^{\sigma_{\lambda}}$ can be found up to quadratures. Let us concisely
see how this works. Suppose for instance that $M$ is a symplectic
manifold with symplectic form $\omega$ and $X$ is a Hamiltonian
vector field defined by a function $H$. Let us write $X=X_{H}$.
In such a case (see Ref. \cite{gp}), the second part of \eqref{Srelr}
is equivalent to 
\begin{equation}
i_{X_{H}^{\Sigma}}\Sigma^{*}\omega=\Sigma^{*}\mathsf{d}H.\label{hjes}
\end{equation}
In other words, $X_{H}^{\Sigma}\in\mathfrak{X}\left(N\times\Lambda\right)$
is a Hamiltonian vector field w.r.t. the symplectic form $\Sigma^{*}\omega$
and with Hamiltonian function $H\circ\Sigma$. Now, the annunciated
condition: $\Sigma$ is said to be \textbf{isotropic} if 
\[
\sigma_{\lambda}^{*}\omega=0,\;\;\;\forall\lambda\in\Lambda.
\]
It was shown in \cite{gp} that, given an isotropic complete solution,
Eq. \eqref{hjes} implies that (unless locally) 
\[
\mathsf{L}_{X_{H}^{\Sigma}}\left(\mathsf{d}W-\theta\right)=\mathfrak{p}_{\Lambda}^{*}\mathsf{d}h,
\]
where $W$ and $h$ are functions related to $\Sigma$, $H$ and $\omega$,
which can be found up to quadratures, and $-\theta$ is a (local)
primitive of $\omega$ (which, by the Poincaré lemma, can also be
found up to quadratures). So, defining for each $\lambda\in\Lambda$
the function $\varphi_{\lambda}:N\rightarrow T_{\lambda}^{*}\Lambda$
such that 
\begin{equation}
\left\langle \varphi_{\lambda}\left(p\right),z\right\rangle =\left\langle \left(\mathsf{d}W-\theta\right)\left(p,\lambda\right),\left(0,z\right)\right\rangle ,\;\;\;\forall z\in T_{\lambda}^{*}\Lambda,\label{dfi}
\end{equation}
we have that each curve $\gamma$ is given by 
\[
\frac{d}{dt}\varphi_{\lambda}\left(\gamma\left(t\right)\right)=\mathsf{d}h\left(\lambda\right),
\]
or equivalently by 
\begin{equation}
\varphi_{\lambda}\left(\gamma\left(t\right)\right)=\varphi_{\lambda}\left(\gamma\left(0\right)\right)+t\,\mathsf{d}h\left(\lambda\right).\label{alfi}
\end{equation}
Since each $\varphi_{\lambda}$ is an immersion, we can solve the
last algebraic equation for $\gamma\left(t\right)$, obtaining in
this way the integrals curves of $X_{H}^{\sigma_{\lambda}}$.

\bigskip{}

In this paper we shall study the above mentioned problem in the case
in which $M$ is a contact manifold and $X$ is a contact vector field.
To do that, we could consider the symplectification of the contact
manifold and then try to apply the procedure described above. But
we wanted to develop our study completely inside the contact manifold
scenario, and this is what we shall do.

The organization of the paper is as follows. In Section 2 we make
a brief review on contact manifolds and contact (Hamiltonian) systems.
In Section 3 we study the partial and complete solutions of the HJE
in the context of contact systems. In Section 4 we introduce the concept
of pseudo-isotropic complete solutions and show that the knowledge
of one of such solutions ensures that the trajectories of the system
can be found up to quadratures. At the end of the section, some comments
are made relating the non-commutative integrability of contact systems,
or \textit{contact non-commutative integrability} (CNCI), and our
pseudo-isotropic complete solutions. In particular, we show that some
of the conditions imposed by the CNCI are not needed to ensure integrability
by quadratures. A similar result was obtained in \cite{gp}, but in
the context of symplectic and Poisson manifolds. The integration procedure
which we use in Section 4, in some sense, works better when the Hamiltonian
function $H$ defining the contact system is a non-vanishing function.
For vanishing Hamiltonian functions, we develop in Section 5 an alternative
procedure. This allowed us to deal with a case that is not explored
(as far as we know) in the literature of integrable contact systems:
the one for which $\xi\left(H\right)\neq0$, where $\xi$ is the Reeb
vector field of the contact manifold.

\bigskip{}

\noindent We assume that the reader is familiar with the basic concepts
of Differential Geometry (see \cite{boot,kn,mrgm}), and in particular
Symplectic Geometry, and with the basic ideas related to Hamiltonian
systems in the context of Geometric Mechanics (see \cite{am,ar,mr}).
We shall work in the smooth (i.e. $C^{\infty}$) category, focusing
exclusively on finite-dimensional smooth manifolds.

\section{Contact Hamiltonian systems}

In this section we recall the definition of contact manifolds, Legendrian
submanifolds, Darboux coordinates for contact manifolds, Reeb vector
fields, contact vector fields, and we say what we shall mean by a
contact Hamiltonian system. For more details, see for instance Refs.
\cite{ar,lib,mru1}.

\subsection{Contact manifolds and contact Hamiltonian vector fields}

By a \textbf{contact manifold} we shall understand a pair $(M,\eta)$,
where $M$ is a $\left(2n+1\right)$-dimensional manifold and $\eta$
is a $1$-form on $M$ such that $\left(\mathsf{d}\eta\right)^{n}\wedge\eta$
is a volume form, i.e. 
\begin{equation}
\begin{array}{ccc}
\underbrace{\textrm{n-times}}\\
\mathsf{(d}\eta)_{p}\wedge\dots\wedge\mathsf{(d}\eta)_{p} & \kern-8pt \wedge & \kern-8pt \eta_{p}\neq0,\;\;\;\forall p\in M.
\end{array}\label{ccm}
\end{equation}
In particular, since $\mathsf{d}\eta\wedge\eta\neq0$, the distribution
$\mathsf{Ker}\,\eta$ is not completely integrable. Moreover, it can
be shown that its maximal integral submanifolds are $n$-dimensional.
Let us say that all the results which we will describe in this section
are proved in Ref. \cite{ar}.

For later convenience, consider the next definition. 
\begin{defn}
A submanifold $S\subset M$ is said to be \textbf{isotropic} (resp.
\textbf{pseudo-isotropic}) if $\mathfrak{i}^{*}\eta=0$ (resp. $\mathfrak{i}^{*}\mathsf{d}\eta=0$),
being $\mathfrak{i}:S\rightarrow M$ the inclusion map. $S$ is said
to be a \textbf{Legendrian} \textbf{submanifold} if is isotropic and
$n$-dimensional. 
\end{defn}
Note that, if $S$ is isotropic, then $\mathsf{Im}\,\mathfrak{i}_{*}=TS\subset\mathsf{Ker}\,\eta$,
i.e. $S$ is an integral submanifold of $\mathsf{Ker}\,\eta$. So,
the Legendrian submanifolds are the maximal integral submanifolds
of $\mathsf{Ker}\,\eta$.

\bigskip{}

Around every $p\in M$ of a contact manifold $M$, there exists a
coordinate chart $\left(U,\left(x^{1},...,x^{n},y_{1},...,y_{n},z\right)\right)$,
called \textbf{Darboux} or \textbf{(contact) canonical coordinate}
\textbf{chart}, such that the local expression of $\eta$ is (sum
over repeated index convention is assumed from now on) 
\[
\eta_{\left|U\right.}=y_{i}\,\mathsf{d}x^{i}+\mathsf{d}z.
\]
In these coordinates, it can be shown that the Legendrian submanifolds
are locally given by equations of the form 
\begin{equation}
y_{i}=-\frac{\partial\Phi}{\partial x^{i}},\;\;\;x^{j}=\frac{\partial\Phi}{\partial y_{j}}\;\;\;\textrm{and}\;\;\;z=\Phi-y_{j}\,\frac{\partial\Phi}{\partial y_{j}},\;\;\;i\in I,\;\;\;j\in J,\label{dar}
\end{equation}
where $I$ and $J$ give a partition of $\left\{ 1,...,n\right\} $
and $\Phi$ is any function of the variables $x^{i}$ and $y_{j}$,
with $i\in I$ and $j\in J$. 
\begin{rem}
\label{canco} Given a cotangent bundle $T^{*}Q$, its canonical $1$-form
$\theta$ induces on $T^{*}Q\times\mathbb{R}$ a contact form given
by $\eta_{Q}=p_{1}^{*}\theta+\mathsf{d}z$, being $p_{1}$ the projection
of $T^{*}Q\times\mathbb{R}$ onto $T^{*}Q$ and $z$ the global coordinate
of the factor $\mathbb{R}$. It is clear that, if $\left(x^{1},...,x^{n},y_{1},...,y_{n}\right)=\left(\mathbf{x},\mathbf{y}\right)$
is a system of canonical coordinates for the symplectic manifold $\left(T^{*}Q,-\mathsf{d}\theta\right)$,
then $\left(\mathbf{x},\mathbf{y},z\right)$ is a system of canonical
coordinates for the contact manifold $\left(T^{*}Q\times\mathbb{R},\eta_{Q}\right)$.
Locally, every contact manifold looks like the last one. 
\end{rem}
Eq. \eqref{ccm} is the same as saying that $\mathsf{Ker}\,\eta$
and $\mathsf{Ker}\,\mathsf{d}\eta$ are regular distributions satisfying
$TM=\mathsf{Ker}\,\eta\oplus\mathsf{Ker}\,\mathsf{d}\eta$. This implies
that $\mathsf{dim}\mathsf{Ker}\,\left(\mathsf{d}\eta\right)_{p}=1$
and $\mathsf{d}\eta_{p}$ is (symplectic) non-degenerated when restricted
to $\mathsf{Ker}\,\eta_{p}$, for all $p\in M.$ As a consequence,
there exists a unique global vector field $\xi\in\mathfrak{X}\left(M\right)$,
called \textbf{Reeb vector field}, characterized by 
\begin{equation}
i_{\xi}\mathsf{d}\eta=0\;\;\;\textrm{and}\;\;\;i_{\xi}\eta=1.\label{drvf}
\end{equation}
Note that $\xi$ is a global generator of $\mathsf{Ker}\,\mathsf{d}\eta$.
Therefore, $TM=\mathsf{Ker}\,\eta\oplus\left\langle \xi\right\rangle $.
Also, in a Darboux chart $U$, we have that the local expression of
the Reeb vector field of $(M,\eta)$ is 
\begin{equation}
\xi_{\left|U\right.}=\frac{\partial}{\partial z}.\label{rd}
\end{equation}
Moreover, given $F\in C^{\infty}\left(M\right)$ and $\alpha\in\Omega^{1}\left(M\right)$,
there exists a unique $X\in\mathfrak{X}\left(M\right)$ such that
\begin{equation}
i_{X}\mathsf{d}\eta=\alpha-i_{\xi}\alpha\,\eta\;\;\;\textrm{and}\;\;\;i_{X}\eta=F.\label{be}
\end{equation}
In particular, if for some function $H\in C^{\infty}\left(M\right)$
we take $F=H$ and $\alpha=-\mathsf{d}H$, Eq. \eqref{be} defines
the so-called \textbf{contact field} $X_{H}$ related to $H$. Thus,
$X_{H}$ is characterized by 
\begin{equation}
i_{X_{H}}\mathsf{d}\eta=-\mathsf{d}H+\xi\left(H\right)\,\eta\;\;\;\textrm{and}\;\;\;i_{X_{H}}\eta=H.\label{dcvf}
\end{equation}
In Darboux coordinates, the local expression of $X_{H}$ is, around
an open neighborhood $U$, 
\begin{equation}
\left(X_{H}\right)_{\left|U\right.}=\left[\frac{\partial H}{\partial y_{i}}\,\frac{\partial}{\partial x^{i}}+\left(y_{i}\,\frac{\partial H}{\partial z}-\frac{\partial H}{\partial x^{i}}\right)\,\frac{\partial}{\partial y_{i}}+\left(H-y_{i}\,\frac{\partial H}{\partial y_{i}}\right)\,\frac{\partial}{\partial z}\right].\label{cfd}
\end{equation}
Note that the Reeb vector field is the contact field related to the
constant function equal to $1$, i.e. $\xi=X_{1}$.

\bigskip{}

Finally, given a contact manifold $\left(M,\eta\right)$ and a function
$H\in C^{\infty}\left(M\right)$, we shall call \textbf{contact Hamiltonian
system}, or simply \textbf{contact system}, to the pair $\left(M,X_{H}\right)$.

\subsection{Some examples}

In this section we present some physical systems that can be described
in terms of contact systems. For more examples, see Ref. \cite{bra}.

\subsubsection{Thermodynamic processes}

\label{tp} Let us consider the contact manifold of Remark \ref{canco}
for $Q=\mathbb{R}^{n}$, i.e. take $M=T^{*}\mathbb{R}^{n}\times\mathbb{R}$.
Let us identify $M$ with $\mathbb{R}^{n}\times\mathbb{R}^{n}\times\mathbb{R}=\mathbb{R}^{2n+1}$.
If $\left(x^{1},...,x^{n},y_{1},...,y_{n},z\right)=\left(\mathbf{x},\mathbf{y},z\right)$
are the global coordinates for $M$, then $\eta_{\mathbb{R}^{n}}=y_{i}\,\mathsf{d}x^{i}+\mathsf{d}z$.
Let us define $H:\mathbb{R}^{2n+1}\rightarrow\mathbb{R}$ as 
\begin{equation}
H\left(\mathbf{x},\mathbf{y},z\right)=a^{0}\left(\mathbf{x}\right)\,\left(z-\Phi\left(\mathbf{x}\right)\right)+a^{j}\left(\mathbf{x}\right)\,\left(y_{j}+\frac{\partial\Phi}{\partial x^{j}}\left(\mathbf{x}\right)\right),\label{Htp}
\end{equation}
for certain functions $a^{j},\Phi:\mathbb{R}^{n}\rightarrow\mathbb{R}$
with some non-null $a^{j}.$ According to \eqref{cfd}, its related
contact field is 
\begin{equation}
\begin{array}{lll}
X_{H} & = & a^{i}\,{\displaystyle \frac{\partial}{\partial x^{i}}+\left[a^{0}\,\left(y_{i}+{\displaystyle \frac{\partial\Phi}{\partial x^{i}}}\right)-{\displaystyle \frac{\partial a^{0}}{\partial x^{i}}\,\left(z-\Phi\right)-\left({\displaystyle \frac{\partial a^{j}}{\partial x^{i}}\,\left(y_{j}+\frac{\partial\Phi}{\partial x^{j}}\right)+a^{j}\,\frac{\partial^{2}\Phi}{\partial x^{i}\partial x^{j}}}\right)}\right]\,\frac{\partial}{\partial y_{i}}}\\
\\
 &  & +\left(a^{0}\,\left(z-\Phi\right)+a^{j}\,{\displaystyle \frac{\partial\Phi}{\partial x^{j}}}\right)\,{\displaystyle \frac{\partial}{\partial z}.}
\end{array}\label{XH}
\end{equation}
It is clear that $\Phi$ defines a Legendrian submanifold $\mathcal{L}\subseteq H^{-1}\left(0\right)$
by the equations {[}recall Eq. \eqref{dar}{]} 
\begin{equation}
y_{i}=-\frac{\partial\Phi}{\partial x^{i}},\;\;\;i=1,...,n,\;\;\;\textrm{and}\;\;\;z=\Phi.\label{L}
\end{equation}
It can be shown that, in general, given a Legendrian submanifold $\mathcal{L}$,
a contact vector field $X_{H}$ is tangent to $\mathcal{L}$ if and
only if $\mathcal{L}\subseteq H^{-1}\left(0\right)$ (see \cite{mru1},
Theorem 3). Then, we have in our case that $\mathsf{Im}\left(X_{H\left|\mathcal{L}\right.}\right)\subseteq T\mathcal{L}$,
for $\mathcal{L}$ given by Eq. \eqref{L}. In addition, taking $\left(x^{1},...,x^{n}\right)$
as coordinates for $\mathcal{L}$, it can be shown that the restriction
of $X_{H}$ to $\mathcal{L}$ is 
\begin{equation}
\hat{X}\left(\mathbf{x}\right)\coloneqq X_{H\left|\mathcal{L}\right.}\left(\mathbf{x}\right)=a^{i}\left(\mathbf{x}\right)\,\left.\frac{\partial}{\partial x^{i}}\right|_{\mathbf{x}}.\label{X}
\end{equation}

\bigskip{}

The physical meaning of $X_{H}$ is as follows (see Ref. \cite{mru1}).
The quasi-static processes of a thermodynamic system are curves on
a Legendrian submanifold $\mathcal{L}$ of some contact manifold.
In many cases, only certain processes are taking into account, given
by the integral curves of some vector field along $\mathcal{L}$.
So, if we have a thermodynamic system with processes contained in
the Legendrian submanifold $\mathcal{L}$ defined by \eqref{L} and
given by the integral curves of the vector field $\hat{X}$ defined
by \eqref{X}, then the integral curves of the contact field $X_{H}$
along $\mathcal{L}$ are precisely the quasi-static processes of the
system under consideration. 
\begin{rem}
\label{a0} Note that the term in $H$ containing the function $a^{0}$
{[}see Eq. \eqref{Htp}{]} plays no role in the physical interpretation
of $X_{H}$ we made above. In particular, we can still make such an
interpretation even though $a^{0}=0$. We add such a term just for
later convenience. 
\end{rem}

\subsubsection{Iso-energetic surfaces and Reeb vector flows}

\label{iso} Let $\left(M,\omega\right)$ be a $\left(2n+2\right)$-dimensional
symplectic manifold. If we have a vector field $\Delta$ on $M$ such
that $\mathsf{L}_{\Delta}\omega=\omega$, we say that $\Delta$ is
a \textbf{Liouville vector field} for $\omega$. (Note that $\mathsf{d}i_{\Delta}\omega=\omega$).
In such a case (see Ref. \cite{lib}), the $1$-form $\eta\coloneqq i_{\Delta}\omega$
determines a contact form on any hyper-surface $S$ of $M$ transverse
to $\Delta,$ i.e. any submanifold $S$ such that 
\[
T_{x}M=T_{x}S\oplus\left\langle \Delta\left(x\right)\right\rangle ,\;\;\;\forall x\in S.
\]
For instance, if the transverse hyper-surface $S$ is given by a level
set of a function $H:M\to{\Bbb R}$ (\textit{ipso facto} $\Delta\left(H\right)$
is non-vanishing along $S$), it is easy to show that the Reeb vector
field $\xi$ associated to $\left(S,i_{\Delta}\omega_{\left|S\right.}\right)$
is 
\begin{equation}
\xi=-\frac{1}{\Delta\left(H\right)_{\left|S\right.}}(X_{H})_{\left|S\right.}.\label{rhx}
\end{equation}
Thus, the trajectories of the (symplectic) Hamiltonian system defined
by $H$, along the given level set $S$, are in essence the integral
curves of a Reeb vector field: the \textbf{Reeb vector flow} (up to
a re-parametrization). In particular, if $\xi$ has closed integral
curves, so does $X_{H}$ along $S$ (see Ref. \cite{w} for more details).

\bigskip{}

In order to give a concrete example (see Ref. \cite{lib}), suppose
that $M=\mathbb{R}^{2n+2}$. Denote by $\left(q^{1},...,q^{n+1},p_{1},...,p_{n+1}\right)$
its global coordinates and consider the symplectic form $\omega$
such that last coordinates are canonical. Then, 
\[
\Delta=\frac{1}{2}\,\left(q^{i}\,\frac{\partial}{\partial q_{i}}+p^{i}\,\frac{\partial}{\partial p_{i}}\right)
\]
is a Liouville vector field for $\omega$. Moreover, the $\left(2n+1\right)$-sphere
$\mathbb{S}^{2n+1}$ is a contact manifold with contact form 
\[
\eta=\frac{1}{2}\,\left(q^{i}\,\mathsf{d}p_{i}-p_{i}\,\mathsf{d}q^{i}\right)_{\left|\mathbb{S}^{2n+1}\right.}
\]
and Reeb vector field
\[
\xi=2\,\left(q^{i}\,\frac{\partial}{\partial p_{i}}-p_{i}\,\frac{\partial}{\partial q^{i}}\right)_{\left|\mathbb{S}^{2n+1}\right.}.
\]
Since $\mathbb{S}^{2n+1}$ can be described as the level set $H=\sum_{i=1}^{n+1}{\displaystyle \left(\left(q^{i}\right)^{2}+p_{i}^{2}\right)=1}$,
and $\Delta\left(H\right)=H$, we have from \eqref{rhx} that $\xi=X_{H\left|\mathbb{S}^{2n+1}\right.}$.
In other words, the integral curves of $X_{H}$ along $\mathbb{S}^{2n+1}$
define exactly a Reeb vector flow.

\section{The $\Pi$-HJE for contact Hamiltonian systems}

Consider a contact manifold $\left(M,\eta\right)$, a contact vector
field $X_{H}$ and a fibration $\Pi:M\rightarrow N$.

\subsection{Partial solutions}

According to Eq. \eqref{hjrel}, having a (partial) solution of the
$\Pi$-HJE for $X_{H}$ is the same as having a section $\sigma:N\to M$
of $\Pi$ satisfying 
\begin{equation}
\sigma_{*}\circ X_{H}^{\sigma}=X_{H}\circ\sigma,\label{hjeh}
\end{equation}
with $X_{H}^{\sigma}\coloneqq\Pi_{*}\circ X_{H}\circ\sigma$. 
\begin{rem}
\label{DS} In Ref. \cite{ds}, a HJE is defined in the restricted
context of the contact manifolds described in Remark \ref{canco}
(the contact form is the same, up to a sign). Such an equation corresponds
exactly to the $\Pi$-HJE with $\Pi=\pi_{Q}\times\mathsf{id}_{\mathbb{R}}$
(in particular, $N=Q\times\mathbb{R}$), where $\pi_{Q}:T^{*}Q\rightarrow Q$
is the canonical cotangent projection. 
\end{rem}
In the following, we shall study Eq. \eqref{hjeh} in some particular
situations. 
\begin{thm}
\label{teo1} Let $\sigma:N\rightarrow M$ be a section of $\Pi$.
If $\sigma$ is a solution of the $\Pi$-HJE for $X_{H}$, then 
\begin{equation}
i_{X_{H}^{\sigma}}\sigma^{*}\mathsf{d}\eta=\sigma^{*}\left(\xi\left(H\right)\,\eta-\mathsf{d}H\right)\;\;\;\textrm{and}\;\;\;i_{X_{H}^{\sigma}}\sigma^{*}\eta=\sigma^{*}H.\label{hjeta}
\end{equation}
On the other hand, if $\sigma$ satisfies above equations, then 
\begin{equation}
\mathsf{Im}\left(\sigma_{*}\circ X_{H}^{\sigma}-X_{H}\circ\sigma\right)\subseteq\mathsf{Ker}\,\Pi_{*}\cap\left(\mathsf{Im}\,\sigma_{*}\right)^{\bot}\cap\mathsf{Ker}\,\eta,\label{imdif}
\end{equation}
where $\bot$ means the orthogonal\footnote{Given a vector space $V$, a subspace $A\subseteq V$ and a (possibly
degenerate) bilinear $\beta:V\times V\rightarrow\mathbb{R}$, the
orthogonal related to $\beta$ is the set $A^{\bot}=\left\{ v\in V:\beta\left(v,a\right)=0,\;\forall a\in A\right\} .$
It is well-known that, for any two subspaces $A,B\subseteq V$, 
\begin{equation}
A^{\bot}\cap B^{\bot}=\left(A+B\right)^{\bot},\;\;\;A^{\bot}+B^{\bot}\subseteq\left(A\cap B\right)^{\bot}\;\;\;\textrm{and}\;\;\;A\subseteq A^{\bot\bot}.\label{wk}
\end{equation}
} w.r.t. the $2$-form $\mathsf{d}\eta$. 
\end{thm}
\textit{Proof.} Given $p\in N$ and a vector $v\in T_{p}N$, 
\begin{equation}
(\mathsf{d}\eta)_{\sigma(p)}\left(\sigma_{*,p}\left(X_{H}^{\sigma}\left(p\right)\right),\sigma_{*,p}\left(v\right)\right)=(\sigma^{*}\mathsf{d}\eta)_{p}\left(X_{H}^{\sigma}\left(p\right),v\right)=\left(i_{X_{H}^{\sigma}}(\sigma^{*}\mathsf{d}\eta)\right)_{p}\left(v\right)\label{p11}
\end{equation}
and, using the first part of \eqref{dcvf}, 
\begin{equation}
(\mathsf{d}\eta)_{\sigma(p)}\left(X_{H}\left(\sigma\left(p\right)\right),\sigma_{*,p}\left(v\right)\right)=(\sigma^{*}\left(i_{X_{H}}\mathsf{d}\eta)\right)_{p}\left(v\right)=(\sigma^{*}\left(\xi\left(H\right)\,\eta-\mathsf{d}H\right))_{p}\left(v\right).\label{p12}
\end{equation}
Thus, if Eq. \eqref{hjeh} holds, the first part of \eqref{hjeta}
immediately follows from above equations. On the other hand, since
\begin{equation}
\eta_{\sigma(p)}\left(\sigma_{*,p}\left(X_{H}^{\sigma}\left(p\right)\right)\right)=\left(i_{X_{H}^{\sigma}}\sigma^{*}\eta\right)_{p}\label{p21}
\end{equation}
and, using the second part of \eqref{dcvf}, 
\begin{equation}
\eta_{\sigma(p)}\left(X_{H}\left(\sigma\left(p\right)\right)\right)=\left(i_{X_{H}}\eta\right)\left(\sigma\left(p\right)\right)=H\left(\sigma\left(p\right)\right)=\sigma^{*}H\left(p\right),\label{p22}
\end{equation}
then the second part of \eqref{hjeta} is obtained if Eq. \eqref{hjeh}
holds.

Let us show the other affirmation of the theorem. First note that,
since $\sigma$ is a section of $\Pi$, then $\Pi_{*}\circ\sigma_{*}=\mathsf{id}_{TN}$
and consequently, using the definition of $X_{H}^{\sigma}$, 
\[
\Pi_{*}\circ\left(\sigma_{*}\circ X_{H}^{\sigma}-X_{H}\circ\sigma\right)=0.
\]
This means that $\mathsf{Im}\left(\sigma_{*}\circ X_{H}^{\sigma}-X_{H}\circ\sigma\right)\subseteq\mathsf{Ker}\,\Pi_{*}$.
Also, subtracting \eqref{p11} and \eqref{p12}, we have from the
first part of \eqref{hjeta} that 
\[
(\mathsf{d}\eta)_{\sigma(p)}\left(\sigma_{*,p}\left(X_{H}^{\sigma}\left(p\right)\right)-X_{H}\left(\sigma\left(p\right)\right),\sigma_{*,p}\left(v\right)\right)=0,
\]
what implies that $\mathsf{Im}\left(\sigma_{*}\circ X_{H}^{\sigma}-X_{H}\circ\sigma\right)\subseteq\left(\mathsf{Im}\,\sigma_{*}\right)^{\bot}$.
Finally, subtracting \eqref{p21} and \eqref{p22}, from the second
part of \eqref{hjeta} we have that 
\[
\eta_{\sigma(p)}\left(\sigma_{*,p}\left(X_{H}^{\sigma}\left(p\right)\right)-X_{H}\left(\sigma\left(p\right)\right)\right)=0,
\]
which ends our proof.$\;\;\;\triangle$

\bigskip{}
 It is worth mentioning that, in general, Eqs. \eqref{hjeh} and \eqref{hjeta}
are not equivalent. Some additional conditions must be fulfilled to
have an equivalence. The following corollaries contain some of those
conditions. 
\begin{cor}
\label{coro1} Eqs. \eqref{hjeh} and \eqref{hjeta} are equivalent
if\footnote{This condition is satisfied by the submersion $\Pi$ used in Ref.
\cite{ds} (see Remark \ref{DS}). } 
\begin{equation}
\mathsf{Ker}\,\Pi_{*}\cap\mathsf{Ker}\,\eta\subseteq\left(\mathsf{Ker}\,\Pi_{*}\right)^{\bot}.\label{qis}
\end{equation}
\end{cor}
\textit{Proof.} It is clear from the last theorem that we only must
prove that \eqref{hjeta} implies \eqref{hjeh}. Also from the last
theorem, we know that \eqref{hjeta} implies \eqref{imdif}. So, it
is enough to prove that \eqref{imdif} implies \eqref{hjeh}. The
hypothesis of the present corollary says that {[}see Eq. \eqref{wk}{]}
\[
\mathsf{Ker}\,\Pi_{*}\cap\left(\mathsf{Im}\,\sigma_{*}\right)^{\bot}\cap\mathsf{Ker}\,\eta\subseteq\left(\mathsf{Ker}\,\Pi_{*}\right)^{\bot}\cap\left(\mathsf{Im}\,\sigma_{*}\right)^{\bot}=\left(\left(\mathsf{Ker}\,\Pi_{*}\right)_{\left|\mathsf{Im}\sigma\right.}+\mathsf{Im}\,\sigma_{*}\right)^{\bot}.
\]
Since $\sigma$ is a section of $\Pi$, we know that 
\begin{equation}
TM_{\left|\mathsf{Im}\sigma\right.}=\left(\mathsf{Ker}\,\Pi_{*}\right)_{\left|\mathsf{Im}\sigma\right.}\oplus\mathsf{Im}\,\sigma_{*},\label{decom}
\end{equation}
and consequently $\mathsf{Ker}\,\Pi_{*}\cap\left(\mathsf{Im}\,\sigma_{*}\right)^{\bot}\cap\mathsf{Ker}\,\eta=\left\{ 0\right\} $.
This fact combined with Eq. \eqref{imdif} drives us precisely to
Eq. \eqref{hjeh}.$\;\;\;\triangle$

\bigskip{}

\begin{cor}
\label{coro2} Eqs. \eqref{hjeh} and \eqref{hjeta} are equivalent
if 
\begin{equation}
\left(\mathsf{Im}\,\sigma_{*}\right)^{\bot}\cap\mathsf{Ker}\,\eta\subseteq\mathsf{Im}\,\sigma_{*}.\label{qcois}
\end{equation}
\end{cor}
\textit{Proof.} The hypothesis says that {[}recall Eq. \eqref{wk}{]}
\[
\mathsf{Ker}\,\Pi_{*}\cap\left(\mathsf{Im}\,\sigma_{*}\right)^{\bot}\cap\mathsf{Ker}\,\eta\subseteq\mathsf{Ker}\,\Pi_{*}\cap\mathsf{Im}\,\sigma_{*}\subseteq\left(\left(\mathsf{Ker}\,\Pi_{*}\right)_{\left|\mathsf{Im}\sigma\right.}^{\bot}+\left(\mathsf{Im}\,\sigma_{*}\right)^{\bot}\right)^{\bot}.
\]
So, our result easily follows by using \eqref{decom} and the same
reasoning as in the proof of the previous corollary. $\;\;\;\triangle$

\bigskip{}

Conditions \eqref{qis} and \eqref{qcois} are similar to the isotropy
and co-isotropy conditions, respectively, of symplectic geometry.
In fact, condition \eqref{qis} is usually called \textbf{pre-isotropy},
and implies the usual isotropy condition for the subbundle $\mathsf{Ker}\,\Pi_{*}\cap\mathsf{Ker}\,\eta\subseteq\mathsf{Ker}\,\eta$
w.r.t. the symplectic fibration on $\mathsf{Ker}\,\eta$ given by
$\mathsf{d}\eta_{\left|\mathsf{Ker}\,\eta\right.}$. 
\begin{defn}
\label{leg} We shall say that a section $\sigma:N\rightarrow M$
is \textbf{isotropic }(resp.\textbf{ Legendrian}) if so is its image
$\mathsf{Im}\sigma$, i.e. if $\mathsf{Im}\,\sigma_{*}\subseteq\mathsf{Ker}\,\eta$,
or equivalently $\sigma^{*}\eta=0$ (resp. and if in addition $2\,\mathsf{dim}N+1=\mathsf{dim}M$). 
\end{defn}
\begin{cor}
If $\sigma$ is an isotropic solution of the $\Pi$-HJE for $X_{H}$,
then $H\circ\sigma=0$. 
\end{cor}
\textit{Proof.} If $\sigma$ is a solution of the $\Pi$-HJE for $X_{H}$,
we know from Theorem \ref{teo1} that \eqref{hjeta} holds. Then,
if $\sigma^{*}\eta=0$, the corollary follows from the second part
of \eqref{hjeta}.$\;\;\;\triangle$

\bigskip{}
 As a direct consequence of the above corollaries we have the next
result. 
\begin{cor}
\label{c4} If some of the conditions \eqref{qis} or \eqref{qcois}
is satisfied and in addition $\sigma$ is isotropic, then Eq. \eqref{hjeh}
as well as Eq. \eqref{hjeta} are equivalent to the equation $H\circ\sigma=0$. 
\end{cor}
\bigskip{}

Let us describe a particular, but very important case. Given a manifold
$Q$, consider the contact manifold $\left(T^{*}Q\times\mathbb{R},\eta_{Q}\right)$
defined in Remark \ref{canco}. Take $\Pi:T^{*}Q\times\mathbb{R}\rightarrow Q$
such that $\Pi=\pi_{Q}\circ p_{1}$ (recall that $\pi_{Q}:T^{*}Q\rightarrow Q$
is the canonical cotangent projection). It is easy to see that $\Pi$
is a fibration and $\mathsf{Ker}\,\Pi_{*}=\left(\mathsf{Ker}\,\Pi_{*}\right)^{\bot}$.
In particular, we have from this last identity that Eq. \eqref{qis}
holds. Let $\sigma:Q\rightarrow T^{*}Q\times\mathbb{R}$ be a section
of $\Pi$. In the canonical coordinates $\left(\mathbf{x},\mathbf{y},z\right)$
(see Remark \ref{canco} again), the local expression of $\sigma$
is given by the formula 
\[
\sigma\left(\mathbf{x}\right)=\left(\mathbf{x},\mathbf{F}\left(\mathbf{x}\right),\Phi\left(\mathbf{x}\right)\right)
\]
for some $\mathbb{R}^{n}$-valued (resp. $\mathbb{R}$-valued) function
$\mathbf{F}$ (resp. $\Phi$). Suppose that $\sigma$ is isotropic
(see Definition \ref{leg}), and consequently Legendrian (for dimensional
reasons). The condition $\sigma^{*}\eta_{Q}=0$ means that 
\[
F_{i}\left(\mathbf{x}\right)+\frac{\partial\Phi}{\partial x^{i}}\left(\mathbf{x}\right)=0.
\]
Since the hypothesis of Corollary \ref{c4} are fulfilled in this
case, the $\Pi$-HJE for any contact field $X_{H}$ reduces to $H\circ\sigma=0$.
In local terms, this means that 
\[
H\left(\mathbf{x},-\nabla\Phi\left(\mathbf{x}\right),\Phi\left(\mathbf{x}\right)\right)=0.
\]
So, restricted to Legendrian sections of $\Pi$, the $\Pi$-HJE represents
a general first order PDE (partial differential equation) defined
by the function $H$.

\bigskip{}

It is well-known that the characteristic field of above PDE is precisely
the contact field $X_{H}$. In particular, if we know the trajectories
of $X_{H}$, then we can construct the solutions of such an equation
(see Ref. \cite{evans}). But we are actually interested in the opposite
problem (not only in the present example, but also in the general
case): finding the trajectories of $X_{H}$ from the solutions of
its $\Pi$-HJE. To do that, it is not enough to have only one solution,
but we need a ``big family'' of them, namely a \textbf{complete
solution} (see the Section 1). Moreover, we need a particular kind
of complete solutions, which we shall present in Section \ref{integ}.

\subsection{Complete solutions}

In this section we study the complete solutions in the case of contact
systems. Let us begin with an example.

\subsubsection{Solutions related to thermodynamic processes}

\label{tps} Let us go back to Section \ref{tp}. Consider the submersion
$\Pi:\mathbb{R}^{2n+1}\rightarrow\mathbb{R}^{n}$ given by $\Pi\left(\mathbf{x},\mathbf{y},z\right)=\mathbf{x}$
and the related $\Pi$-HJE for $X_{H}$ {[}see Eq. \eqref{hjeh}{]}.
It is easy to show that 
\begin{equation}
X_{H}^{\sigma}=a^{i}\,\frac{\partial}{\partial x^{i}}\label{inds}
\end{equation}
for every section $\sigma:\mathbb{R}^{n}\rightarrow\mathbb{R}^{2n+1}$
of $\Pi$. If we write 
\[
\sigma\left(\mathbf{x}\right)=\left(\mathbf{x},\sigma^{y}\left(\mathbf{x}\right),\sigma^{z}\left(\mathbf{x}\right)\right),
\]
then Eq. \eqref{hjeh} translates to 
\begin{equation}
a^{j}\,\frac{\partial\sigma_{i}^{y}}{\partial x^{j}}=a^{0}\,\left(\sigma_{i}^{y}+\frac{\partial\Phi}{\partial x^{i}}\right)-\frac{\partial a^{0}}{\partial x^{i}}\,\left(\sigma^{z}-\Phi\right)-\left[\frac{\partial a^{j}}{\partial x^{i}}\,\left(\sigma_{j}^{y}+\frac{\partial\Phi}{\partial x^{j}}\right)+a^{j}\,\frac{\partial^{2}\Phi}{\partial x^{i}\partial x^{j}}\right]\label{1}
\end{equation}
and 
\begin{equation}
a^{j}\,\frac{\partial\sigma^{z}}{\partial x^{j}}=a^{0}\,\left(\sigma^{z}-\Phi\right)+a^{j}\,\frac{\partial\Phi}{\partial x^{j}}.\label{2}
\end{equation}
Defining $\alpha_{i}\left(\mathbf{x}\right)\coloneqq\sigma_{i}^{y}\left(\mathbf{x}\right)+{\displaystyle \frac{\partial\Phi}{\partial x^{i}}\left(\mathbf{x}\right)}$
and $\beta\left(\mathbf{x}\right)\coloneqq\sigma^{z}\left(\mathbf{x}\right)-\Phi\left(\mathbf{x}\right)$,
the Eqs. \eqref{1} and \eqref{2} transform into 
\[
a^{j}\,\left(\frac{\partial\alpha_{i}}{\partial x^{j}}-\frac{\partial\alpha_{j}}{\partial x^{i}}\right)=a^{0}\,\alpha_{i}-\frac{\partial a^{0}}{\partial x^{i}}\,\beta-\frac{\partial}{\partial x^{i}}\left(a^{j}\,\alpha_{j}\right)
\]
and 
\[
a^{j}\,\frac{\partial\beta}{\partial x^{j}}=a^{0}\,\beta,
\]
respectively. If we think of the functions $a_{i}$ and $\alpha_{i}$
as the components of a vector field and a $1$-form 
\begin{equation}
\hat{X}\coloneqq a^{i}\,\frac{\partial}{\partial x^{i}}\;\;\;\textrm{and}\;\;\;\alpha\coloneqq\alpha_{i}\,\mathsf{d}x^{i},\label{xa}
\end{equation}
respectively, we can write above equations as 
\begin{equation}
\mathsf{L}_{\hat{X}}\alpha=a^{0}\,\alpha-\beta\,\mathsf{d}a^{0}\label{eq:1ra}
\end{equation}
and 
\begin{equation}
\mathsf{L}_{\hat{X}}\beta=a^{0}\,\beta.\label{eq:2da}
\end{equation}
Let us suppose that $a^{0}$ is \textbf{constant} and take $\alpha=e^{-f}\,\mathsf{d}g$.
Then, Eq. (\ref{eq:1ra}) reduces to 
\[
\mathsf{L}_{\hat{X}}\left(e^{-f}\right)\,\mathsf{d}g+e^{-f}\,\mathsf{L}_{\hat{X}}\left(\mathsf{d}g\right)=a^{0}\,e^{-f}\,\mathsf{d}g,
\]
or equivalently 
\[
-\mathsf{L}_{\hat{X}}\left(f\right)\,\mathsf{d}g+\mathsf{L}_{\hat{X}}\left(\mathsf{d}g\right)=a^{0}\,\mathsf{d}g.
\]
Choosing $f$ such that 
\begin{equation}
\mathsf{L}_{\hat{X}}\left(f\right)=\mathfrak{i}_{\hat{X}}\mathsf{d}f=-a^{0},\label{lxf}
\end{equation}
the function $g$ must satisfy 
\begin{equation}
\mathsf{L}_{\hat{X}}\left(\mathsf{d}g\right)=\mathsf{d}\left(\mathfrak{i}_{\hat{X}}\mathsf{d}g\right)=0.\label{fg}
\end{equation}

\begin{rem}
\label{a001} Note that, if $a^{0}=0$, then we can take $f=0$. 
\end{rem}
Since $\hat{X}$ is a non-null vector field, it is easy to show that
given $n$ constant $c_{1},\dots,c_{n}$, some $c_{k}$ non-null,
then we can find $n$ independent functions $\left\{ g_{1},...,g_{n}\right\} $
satisfying 
\begin{equation}
\mathfrak{i}_{\hat{X}}\mathsf{d}g_{k}=c_{k},\;\;\;k=1,...,n.\label{ck}
\end{equation}
Hence, we have $n$ independent functions solving Eq. \eqref{fg}. 
\begin{rem}
\label{cm} Suppose that $c_{1}\neq0$. Then the $n-1$ functions
$g_{i}-{\displaystyle \frac{c_{i}}{c_{1}}\,g_{1}}$, $i=2,...,n$,
are first integrals for the vector field $\hat{X}$ {[}see Eq. \eqref{ck}{]}. 
\end{rem}
Finally, given constant $\lambda^{1},...,\lambda^{n}$, it easily
follows that $\alpha\coloneqq e^{-f}\,\lambda^{k}\,\mathsf{d}g_{k}$
solves \eqref{eq:1ra}. Also, given any constant $\lambda^{0}$, the
function $\beta\coloneqq\lambda^{0}\,e^{-f}$ solves \eqref{eq:2da}.
In this way, we have a family of solutions of the $\Pi$-HJE for $X_{H}$
given by 
\begin{equation}
\sigma_{\lambda}\left(\mathbf{x}\right)\coloneqq\left(\mathbf{x},e^{-f\left(\mathbf{x}\right)}\,\left(\lambda^{k}\,\nabla g_{k}\left(\mathbf{x}\right)\right)-\nabla\Phi\left(\mathbf{x}\right),\Phi\left(\mathbf{x}\right)+\lambda^{0}\,e^{-f\left(\mathbf{x}\right)}\right),\label{sil}
\end{equation}
where $\lambda=\left(\lambda^{0},\lambda^{1},...,\lambda^{n}\right)\in\mathbb{R}^{n+1}$.
Such a family defines a complete solution 
\begin{equation}
\Sigma:\mathbb{R}^{n}\times\mathbb{R}^{n+1}\rightarrow\mathbb{R}^{2n+1}:\left(\mathbf{x},\lambda\right)\mapsto\sigma_{\lambda}\left(\mathbf{x}\right).\label{Sil}
\end{equation}
In fact, the independence of the functions $g_{i}$'s ensures that
$\Sigma$ is a surjective local diffeomorphism. 
\begin{rem}
\label{isoa0} According to Remark \ref{a001}, if $a^{0}=0$ we have
a complete solution given by 
\[
\sigma_{\lambda}\left(\mathbf{x}\right)\coloneqq\left(\mathbf{x},\lambda^{k}\,\nabla g_{k}\left(\mathbf{x}\right)-\nabla\Phi\left(\mathbf{x}\right),\Phi\left(\mathbf{x}\right)+\lambda^{0}\right).
\]
\end{rem}
Note that the form of $X_{H}^{\sigma_{\lambda}}$ does not depend
on $\lambda$ {[}see Eq. \eqref{inds}{]}. On the other hand, for
$\lambda=\left(\lambda^{0},\lambda^{1},...,\lambda^{n}\right)=\left(0,0,...,0\right)\equiv0$,
we have that $\mathsf{Im}\sigma_{0}=\mathcal{L}$ {[}see Eq. \eqref{L}{]}
and $X_{H}^{\sigma_{0}}=X_{H\left|\mathcal{L}\right.}=\hat{X}$ {[}see
Eqs. \eqref{X} and \eqref{xa}{]}. So, one of the partial solutions
defines a Legendrian submanifold.
\begin{rem*}
A relationship between thermodynamical systems and the Hamilton-Jacobi
Theory was also explored in Ref. \cite{raj}, but from another point
of view.
\end{rem*}

\subsubsection{A useful characterization}

In the following, we shall give an alternative description of the
complete solutions. Fix a contact manifold $(M,\eta)$ and a fibration
$\Pi:M\rightarrow N$. 
\begin{thm}
\label{t2} A surjective local diffeomorphism $\Sigma:N\times\Lambda\rightarrow M$,
satisfying $\Pi\circ\Sigma=\mathfrak{p}_{N}$, is a complete solution
of the $\Pi$-HJE for a contact vector $X_{H}$ if and only if 
\begin{equation}
i_{X_{H}^{\Sigma}}\Sigma^{*}\mathsf{d}\eta=\Sigma^{*}\left(\xi\left(H\right)\,\eta-\mathsf{d}H\right)\;\;\;\textrm{and}\;\;\;i_{X_{H}^{\Sigma}}\Sigma^{*}\eta=\Sigma^{*}H,\label{ccont}
\end{equation}
i.e. $X_{H}^{\Sigma}$ is a contact field for the contact form $\Sigma^{*}\eta$
and for the function $\Sigma^{*}H=H\circ\Sigma$. 
\end{thm}
\textit{Proof.} Suppose that $\Sigma$ is a complete solution. Since
\begin{equation}
\left(i_{X_{H}^{\Sigma}}\Sigma^{*}\eta\right)\left(p,\lambda\right)=\eta_{\Sigma(p,\lambda)}\left(\Sigma_{*}\circ X_{H}^{\Sigma}\left(p,\lambda\right)\right),\;\;\;(p,\lambda)\in N\times\Lambda,\label{au11}
\end{equation}
the second part of \eqref{Srelr} tells us that 
\[
\left(i_{X_{H}^{\Sigma}}\Sigma^{*}\eta\right)\left(p,\lambda\right)=\eta_{\Sigma(p,\lambda)}\left(X_{H}\circ\Sigma\left(p,\lambda\right)\right)=\left(i_{X_{H}}\eta\right)\left(\Sigma\left(p,\lambda\right)\right).
\]
Then, we obtain the second part of \eqref{ccont} from the corresponding
one of \eqref{dcvf}. On the other hand, given any vector $v\in T_{p}N\times T_{\lambda}\Lambda$,
\begin{equation}
\left(i_{X_{H}^{\Sigma}}\Sigma^{*}\mathsf{d}\eta\right)_{(p,\lambda)}\left(v\right)=(\mathsf{d}\eta)_{\Sigma(p,\lambda)}\left(\Sigma_{*}\circ X_{H}^{\Sigma}\left(p,\lambda\right),\Sigma_{*,\left(p,\lambda\right)}\left(v\right)\right),\label{au1}
\end{equation}
and, using the second part of \eqref{Srelr}, 
\[
\left(i_{X_{H}^{\Sigma}}\Sigma^{*}\mathsf{d}\eta\right)_{(p,\lambda)}\left(v\right)=(\mathsf{d}\eta)_{\Sigma(p,\lambda)}\left(X_{H}\circ\Sigma\left(p,\lambda\right),\Sigma_{*,\left(p,\lambda\right)}\left(v\right)\right)=\left(i_{X_{H}}\mathsf{d}\eta\right)_{\Sigma(p,\lambda)}\left(\Sigma_{*,\left(p,\lambda\right)}\left(v\right)\right).
\]
Also, 
\begin{equation}
\left[\Sigma^{*}\left(\xi\left(H\right)\,\eta-\mathsf{d}H\right)\right]_{(p,\lambda)}\left(v\right)=\left(\xi\left(H\right)\,\eta-\mathsf{d}H\right)_{\Sigma(p,\lambda)}\left(\Sigma_{*,\left(p,\lambda\right)}\left(v\right)\right).\label{au2}
\end{equation}
So, from the first part of \eqref{dcvf} we obtain the first one of
\eqref{ccont}.

Now, let us show the converse. Suppose that a local diffeomorphism
$\Sigma$ satisfies \eqref{ccont}. From \eqref{au11} and the second
parts of \eqref{dcvf} and \eqref{ccont} it follows that 
\begin{equation}
\eta_{\Sigma(p,\lambda)}\left(\Sigma_{*}\circ X_{H}^{\Sigma}\left(p,\lambda\right)-X_{H}\circ\Sigma\left(p,\lambda\right)\right)=0.\label{e0}
\end{equation}
On the other hand, given again $v\in T_{p}N\times T_{\lambda}\Lambda$,
Eqs. \eqref{au1} and \eqref{au2} and the first parts of \eqref{dcvf}
and \eqref{ccont} tell us that 
\[
(\mathsf{d}\eta)_{\Sigma(p,\lambda)}\left(\Sigma_{*}\circ X_{H}^{\Sigma}\left(p,\lambda\right)-X_{H}\circ\Sigma\left(p,\lambda\right),\Sigma_{*,\left(p,\lambda\right)}\left(v\right)\right)=0.
\]
Since $v$ is arbitrary and $\Sigma_{*,\left(p,\lambda\right)}$ is
bijective, above equation implies that there is $\alpha\in{\Bbb R}$
such that 
\begin{equation}
\Sigma_{*}\circ X_{H}^{\Sigma}\left(p,\lambda\right)-X_{H}\circ\Sigma\left(p,\lambda\right)=\alpha\,\xi\left(\Sigma\left(p,\lambda\right)\right).\label{e1}
\end{equation}
Combining \eqref{e0} and \eqref{e1}, we have that $\alpha=0$, from
which the converse of the theorem follows.$\;\;\;\triangle$

\bigskip{}

It is worth noticing that being a complete solution is (in essence)
the same as fulfilling Eq. \eqref{ccont} (as we showed above), but
being a partial solution is not the same as satisfying Eq. \eqref{hjeta}
(see the corollaries of the previous section).

\section{An integration procedure}

\label{integ}Fix again a contact manifold $\left(M,\eta\right)$,
a contact vector field $X_{H}$ and a fibration $\Pi:M\rightarrow N$.
We shall see in this section that, under some additional assumptions,
having a complete solution of the $\Pi$-HJE for $X_{H}$ ensures
that the integral curves of $X_{H}$ can be constructed up to quadratures.
This extends to contact manifolds the results obtained in \cite{gp}
(and briefly reviewed in Section 1) for symplectic and Poisson manifolds.

\subsection{Pseudo-isotropy condition}

Assume that we have a complete solution $\Sigma:N\times\Lambda\rightarrow M$
of the $\Pi$-HJE for some vector field $X$, non-necessarily Hamiltonian. 
\begin{defn}
\label{isotp} We shall say that $\Sigma$ is \textbf{pseudo-isotropic}
if 
\begin{equation}
\sigma_{\lambda}^{*}\mathsf{d}\eta=0,\;\;\;\forall\lambda\in\Lambda,\label{isco}
\end{equation}
and that it is \textbf{isotropic} if each partial solution $\sigma_{\lambda}$
is isotropic, i.e. 
\begin{equation}
\sigma_{\lambda}^{*}\eta=0,\;\;\;\forall\lambda\in\Lambda.\label{iso}
\end{equation}
\end{defn}
Note that isotropy implies pseudo-isotropy. Also, the condition \eqref{isco}
is the same as saying that the partial solutions satisfy $\mathsf{Im}\left(\sigma_{\lambda}\right)_{*}\subseteq\left(\mathsf{Im}\left(\sigma_{\lambda}\right)_{*}\right)^{\bot}$. 
\begin{rem*}
It is worth mentioning that isotropy condition does not hold, generically,
for contact Hamiltonian systems. In fact, if $\sigma_{\lambda}$ were
isotropic for all $\lambda$, the second part of Eq. \eqref{hjeta}
would say that $H$ is the null function. Nevertheless, see Definition
\ref{bis}.
\end{rem*}
Let us go back to the case of contact Hamiltonian systems. An immediate
consequence of Eq. \eqref{isco} is the fact that, for all $\lambda\in\Lambda,$
we can construct up to quadratures a function $W_{\lambda}$ such
that $\sigma_{\lambda}^{*}\eta=\mathsf{d}W_{\lambda}$, at least locally.
From now on, assume for simplicity that $N$ is connected and simply
connected. 
\begin{rem}
\label{csc} If this is not the case, for any $p\in N$ consider a
connected and simply connected open neighborhood $V$ of $p$ and
an open subset $\Lambda_{V}\subseteq\Lambda$ such that $\Sigma$
restricted to $V\times\Lambda_{V}$ is a diffeomorphism. Then, in
what follows, replace $N$ by $V$, $\Lambda$ by $\Lambda_{V}$,
and $M$ by $\Sigma\left(V\times\Lambda_{V}\right)$. 
\end{rem}
Under above assumption, $W_{\lambda}$ is globally defined and the
function $W:N\times\Lambda\rightarrow\mathbb{R}:\left(p,\lambda\right)\rightarrow W_{\lambda}\left(p\right)$
satisfies 
\begin{equation}
i_{X}\left(\mathsf{d}W-\Sigma^{*}\eta\right)=0,\label{mai}
\end{equation}
for all $X\in\mathfrak{X}\left(N\times\Lambda\right)$ such that $\mathsf{Im}X\subseteq TN\times\left\{ 0\right\} $.
In particular, since $\mathsf{Im}X_{H}^{\Sigma}\subseteq TN\times\left\{ 0\right\} $
{[}recall Eq. \eqref{xsin}{]}, we have from \eqref{ccont} and \eqref{mai}
that 
\begin{equation}
\mathsf{L}_{X_{H}^{\Sigma}}W=i_{X_{H}^{\Sigma}}\mathsf{d}W=i_{X_{H}^{\Sigma}}\Sigma^{*}\eta=\Sigma^{*}H\label{ma}
\end{equation}
and 
\begin{equation}
i_{X_{H}^{\Sigma}}\Sigma^{*}\mathsf{d}\eta=i_{X_{H}^{\Sigma}}\mathsf{d}\Sigma^{*}\eta=-i_{X_{H}^{\Sigma}}\mathsf{d}\left(\mathsf{d}W-\Sigma^{*}\eta\right)=-\mathsf{L}_{X_{H}^{\Sigma}}\left(\mathsf{d}W-\Sigma^{*}\eta\right).\label{mi}
\end{equation}
Then, on the one hand, for every $\lambda\in\Lambda$ and every integral
curve $\gamma$ of $X_{H}^{\sigma_{\lambda}}$, it follows from \eqref{ma}
that 
\begin{equation}
\frac{d}{dt}W_{\lambda}\left(\gamma\left(t\right)\right)=H\left(\Sigma\left(\gamma\left(t\right),\lambda\right)\right).\label{ma2}
\end{equation}
On the other hand, consider the function $\varphi:N\times\Lambda\rightarrow T^{*}\Lambda$
such that $\varphi\left(p,\lambda\right)\in T_{\lambda}^{*}\Lambda$,
for all $\left(p,\lambda\right)\in N\times\Lambda$, and 
\begin{equation}
\left\langle \varphi\left(p,\lambda\right),z\right\rangle =\left\langle \left(\mathsf{d}W-\Sigma^{*}\eta\right)\left(p,\lambda\right),\left(0,z\right)\right\rangle ,\;\;\;\forall z\in T_{\lambda}\Lambda.\label{defi}
\end{equation}
(Recall Eq. \eqref{dfi} of Section 1). Consider also the related
functions $\varphi_{\lambda}:p\in N\mapsto\varphi\left(p,\lambda\right)\in T_{\lambda}^{*}\Lambda$.
Then, using \eqref{mai} and the fact that $\mathsf{Im}\left[X_{H}^{\Sigma},Z\right]\subseteq TN\times\left\{ 0\right\} $
for any vector field $Z\in\mathfrak{X}\left(\Lambda\right)$, we have
that
\begin{equation}
\mathsf{L}_{X_{H}^{\Sigma}}\circ i_{Z}\left(\mathsf{d}W-\Sigma^{*}\eta\right)=i_{Z}\circ\mathsf{L}_{X_{H}^{\Sigma}}\left(\mathsf{d}W-\Sigma^{*}\eta\right)\label{conm}
\end{equation}
and consequently {[}see Eq. \eqref{mi}{]} 
\begin{equation}
\left\langle \frac{d}{dt}\varphi_{\lambda}\left(\gamma\left(t\right)\right),z\right\rangle =\left\langle \left(\mathsf{L}_{X_{H}^{\Sigma}}\left(\mathsf{d}W-\Sigma^{*}\eta\right)\right)\left(\gamma\left(t\right),\lambda\right),\left(0,z\right)\right\rangle =-\left\langle \left(i_{X_{H}^{\Sigma}}\Sigma^{*}\mathsf{d}\eta\right)\left(\gamma\left(t\right),\lambda\right),\left(0,z\right)\right\rangle ,\label{mi2}
\end{equation}
for all $z\in T_{\lambda}\Lambda$. Identities \eqref{ma2} and \eqref{mi2},
together with the next result, will be very useful in the rest of
the paper. 
\begin{prop}
\label{inm} Each function$\left(\varphi_{\lambda},W_{\lambda}\right):N\rightarrow T_{\lambda}^{*}\Lambda\times\mathbb{R}$
is an immersion. As a consequence, 
\[
\left(\varphi,W\right):N\times\Lambda\rightarrow T^{*}\Lambda\times\mathbb{R}
\]
 is an immersion also. 
\end{prop}
\textit{Proof}. Since $\varphi_{\lambda}$ takes values in the vector
space $T_{\lambda}^{*}\Lambda$, let us identify its differential
$\left(\varphi_{\lambda}\right)_{*}$ with a function from $TN$ to
$T_{\lambda}^{*}\Lambda$. Under this identification, it is easy to
show that 
\[
\left\langle \left(\varphi_{\lambda}\right)_{*,p}\left(x\right),z\right\rangle =\left[\mathsf{d}\left(\mathsf{d}W-\Sigma^{*}\eta\right)\right]_{\left(p,\lambda\right)}\left(\left(x,0\right),\left(0,z\right)\right)=-\left(\Sigma^{*}\mathsf{d}\eta\right)_{\left(p,\lambda\right)}\left(\left(x,0\right),\left(0,z\right)\right),
\]
for all $p\in N$, $x\in T_{p}N$ and $z\in T_{\lambda}\Lambda$.
Note that, by the pseudo-isotropy condition, 
\[
\left(\Sigma^{*}\mathsf{d}\eta\right)_{\left(p,\lambda\right)}\left(\left(x,0\right),\left(0,z\right)\right)=\left(\Sigma^{*}\mathsf{d}\eta\right)_{\left(p,\lambda\right)}\left(\left(x,0\right),\left(y,z\right)\right),
\]
for all $y\in T_{p}N$. Then, if $\left(\varphi_{\lambda}\right)_{*,p}\left(x\right)=0$,
we must have that $\Sigma_{*,(p,\lambda)}\left(x,0\right)=\alpha\,\xi\left(p,\lambda\right)$,
for some number $\alpha$. On the other hand, also by pseudo-isotropy
{[}see for instance Eq. \eqref{mai}{]}, 
\[
\left(W_{\lambda}\right)_{*,p}\left(x\right)=\left(\Sigma^{*}\eta\right)_{\left(p,\lambda\right)}\left(x,0\right)=\alpha.
\]
So, if we also ask that $\left(W_{\lambda}\right)_{*,p}\left(x\right)=0$,
then $\alpha=0$ and consequently $x$ must be zero. This ends our
proof.$\;\;\;\triangle$

\bigskip{}

The theorem below is the analogous for contact manifolds of Proposition
3.17 of Ref. \cite{gp} (for symplectic manifolds), and reflects the
geometric meaning of the function $\left(\varphi,W\right)$. 
\begin{thm*}
Let $(M,\eta)$ be a contact manifold and $\Sigma:N\times\Lambda\rightarrow M$
a pseudo-isotropic complete solution of the $\Pi$-HJE for some vector
field $X.$ Consider on $N\times\Lambda$ the contact form $\Sigma^{*}\eta$
and on $T^{*}\Lambda\times\mathbb{R}$ the one given by the canonical
$1$-form on $T^{*}\Lambda$ (see Remark \ref{canco}). Then, $\left(\varphi,W\right):N\times\Lambda\rightarrow T^{*}\Lambda\times\mathbb{R}$
is a anti-morphism of contact manifolds.\footnote{Given two contact manifolds $\left(M_{1},\eta_{1}\right)$ and $\left(M_{2},\eta_{2}\right)$,
a smooth function $f:M_{1}\rightarrow M_{2}$ is a \textbf{morphism}
(resp. \textbf{anti-morphism}) of contact manifolds if $f^{*}\eta_{2}=\eta_{1}$
(resp. $f^{*}\eta_{2}=-\eta_{1}$).} 
\end{thm*}
\textit{Proof}. We must show that $\left(\varphi,W\right)^{*}\eta_{\Lambda}=-\Sigma^{*}\eta$.
Given a vector $\left(x,y\right)\in T_{p}N\times T_{\lambda}\Lambda$,
we have that 
\[
\begin{array}{lll}
\left[\left(\varphi,W\right)^{*}\eta_{\Lambda}\right]\left(x,y\right) & = & \eta_{\Lambda}\left(\varphi_{*,\left(p,\lambda\right)}\left(x,y\right),\mathsf{d}W\left(p,\lambda\right)\left(x,y\right)\right)\\
 & = & \left\langle \varphi_{\lambda}\left(p\right),y\right\rangle +\left\langle \mathsf{d}W\left(p,\lambda\right),\left(x,y\right)\right\rangle \\
 & = & \left\langle \left(\mathsf{d}W-\Sigma^{*}\eta\right)\left(p,\lambda\right),\left(0,y\right)\right\rangle +\left\langle \mathsf{d}W\left(p,\lambda\right),\left(x,y\right)\right\rangle ,
\end{array}
\]
and, using that $\left\langle \left(\mathsf{d}W-\Sigma^{*}\eta\right)\left(p,\lambda\right),\left(x,0\right)\right\rangle =0$
{[}see Eq. \eqref{mai}{]}, the wanted identity follows.$\;\;\;\triangle$

\subsection{Contact systems satisfying $\xi\left(H\right)=0$ }

\label{si0} Now suppose that we have a contact system with $H$ such
that $\xi\left(H\right)=0$. Of course, this condition is fulfilled
by the Reeb vector field. 
\begin{rem*}
In the example of Section \ref{tp}, we have that $\xi\left(H\right)=a^{0}$
{[}see Eqs. \eqref{rd} and \eqref{Htp}{]}. So, taking $a^{0}=0$,
we have from \eqref{Htp} a big family of examples satisfying the
condition $\xi\left(H\right)=0$, all of them with a precise physical
meaning (see Remark \ref{a0}). 
\end{rem*}
Let $\Sigma$ be a pseudo-isotropic complete solution. When $\xi\left(H\right)=0$,
Eq. \eqref{ccont} says that $i_{X_{H}^{\Sigma}}\Sigma^{*}\mathsf{d}\eta=-\Sigma^{*}\mathsf{d}H$.
So, by the pseudo-isotropy condition we have that
\[
i_{X}\Sigma^{*}\mathsf{d}H=0,
\]
for all $X\in\mathfrak{X}\left(N\times\Lambda\right)$ such that $\mathsf{Im}X\subseteq TN\times\left\{ 0\right\} $,
or equivalently 
\[
\sigma_{\lambda}^{*}\mathsf{d}H=0,\;\;\;\forall\lambda\in\Lambda.
\]
Assuming again that $N$ is connected, it follows that each function
$H\circ\sigma_{\lambda}\in C^{\infty}\left(N\right)$ is constant.
In other words, we can construct a function $h:\Lambda\rightarrow\mathbb{R}$
such that 
\begin{equation}
h\circ\mathfrak{p}_{\Lambda}=H\circ\Sigma.\label{h}
\end{equation}
This implies that 
\[
i_{X_{H}^{\Sigma}}\Sigma^{*}\mathsf{d}\eta=-\mathfrak{p}_{\Lambda}^{*}\mathsf{d}h.
\]
As a consequence, given an integral curve $\gamma$ of $X_{H}^{\sigma_{\lambda}}$,
Eq. \eqref{mi2} tells us in this case that 
\[
\frac{d}{dt}\varphi_{\lambda}\left(\gamma\left(t\right)\right)=\mathsf{d}h\left(\lambda\right),
\]
or equivalently, 
\begin{equation}
\varphi_{\lambda}\left(\gamma\left(t\right)\right)=t\,\mathsf{d}h\left(\lambda\right)+\varphi_{\lambda}\left(\gamma\left(0\right)\right).\label{fsh}
\end{equation}
The latter is exactly the equation that appeared in Section 1 {[}see
Eq. \eqref{alfi}{]} in the context of symplectic manifolds. On the
other hand, using \eqref{ma2} and \eqref{h}, 
\[
\frac{d}{dt}W_{\lambda}\left(\gamma\left(t\right)\right)=h\left(\lambda\right),
\]
i.e. 
\begin{equation}
W_{\lambda}\left(\gamma\left(t\right)\right)=t\,h\left(\lambda\right)+W_{\lambda}\left(\gamma\left(0\right)\right).\label{Wsh}
\end{equation}
Accordingly, if we want to find the integral curves $\gamma$ of each
vector field $X_{H}^{\sigma_{\lambda}}$, we just need to solve the
algebraic equations given by \eqref{fsh} and \eqref{Wsh}, which
can be solved univocally because the function $\left(\varphi_{\lambda},W_{\lambda}\right)$
is an immersion (see Proposition \ref{inm}). Summing up, we have
the next result. 
\begin{thm}
\label{s0pi}Given a contact manifold $\left(M,\eta\right)$, a fibration
$\Pi:M\rightarrow N$ and a function $H:M\rightarrow\mathbb{R}$ such
that $\xi\left(H\right)=0$, if we know a pseudo-isotropic complete
solution $\Sigma$ of the $\Pi$-HJE for $X_{H}$, then the integral
curves of $X_{H}$ can be constructed up to quadratures. In the particular
case when $H=1$, we have that the integral curves of $X_{H}=\xi$,
i.e. the Reeb vector flow, can be constructed up to quadratures. 
\end{thm}
Its proof follows directly from above calculations when $N$ is connected
and simply connected. Otherwise, we can replicate such calculations
around every point $p\in N$, as emphasized in Remark \ref{csc}.\bigskip{}

The previous theorem extends to contact systems $\left(M,X_{H}\right)$
(with $\xi\left(H\right)=0$) the result obtained in Ref. \cite{gp}
for Hamiltonian systems on symplectic manifolds. When the condition
$\xi\left(H\right)=0$ is not satisfied, we need another kind of complete
solutions. They will be studied in the next subsection.

\subsection{The $g$-pseudo-isotropy condition}

\label{ges} Consider again a complete solution $\Sigma:N\times\Lambda\rightarrow M$
of the $\Pi$-HJE for some contact field $X_{H}$. 
\begin{defn}
\label{geta} We shall say that $\Sigma$ is $g$-\textbf{pseudo-isotropic
for $H$} if 
\begin{equation}
\sigma_{\lambda}^{*}\mathsf{d}\left(g\,\eta\right)=0,\;\;\;\forall\lambda\in\Lambda,\label{isco-1}
\end{equation}
for some non-vanishing function $g\in C^{\infty}\left(M\right)$ such
that 
\begin{equation}
X_{H}\left(g\right)+g\,\xi\left(H\right)=0.\label{lxhg}
\end{equation}
\end{defn}
A pseudo-isotropic solution is also $g$-pseudo-isotropic for every
$H$ such that $\xi\left(H\right)=0$. We just must take $g=1$ above.

The next result try to explain the meaning of Eq. \eqref{lxhg}. The
proof is left to the reader. 
\begin{prop}
\label{modcon} Given two functions $H$ and $g$ on $M$, we have
that 
\begin{equation}
i_{X_{H}}\mathsf{d}\left(g\eta\right)=\left(X_{H}\left(g\right)+g\,\xi\left(H\right)\right)\,\eta-\mathsf{d}\left(gH\right)\;\;\;\textrm{and}\;\;\;i_{X_{H}}g\eta=gH.\label{gheta}
\end{equation}
In addition, if $g$ is a non-vanishing function, then $X_{H}$ is
also a contact field w.r.t. to the contact form $\eta'=g\eta$, but
with Hamiltonian function $H'=gH$. And denoting $\xi'$ the Reeb
vector field for $\eta'$, we have that 
\begin{equation}
\xi'\left(H'\right)=\frac{X_{H}\left(g\right)+g\,\xi\left(H\right)}{g}.\label{sphp}
\end{equation}
\end{prop}
So, if we have a $g$-pseudo-isotropic complete solution for $H$,
we can see it as a pseudo-isotropic complete solution for the contact
manifold $\left(M,\eta'\right)$ and the contact field with Hamiltonian
$H'$. Moreover, from Eqs. \eqref{lxhg} and \eqref{sphp} we have
that $\xi'\left(H'\right)=0$. Consequently, we are in the situation
of the Theorem \ref{s0pi}, and the next result easily follows. 
\begin{thm}
\label{impo} Given a contact manifold $\left(M,\eta\right)$, a fibration
$\Pi:M\rightarrow N$ and a function $H:M\rightarrow\mathbb{R}$,
if we know a $g$-pseudo-isotropic complete solution $\Sigma$ of
the $\Pi$-HJE for $X_{H}$, then the integral curves of $X_{H}$
can be constructed up to quadratures. 
\end{thm}
\bigskip{}

Suppose now that we have a contact vector field $X_{H}$ with $H\left(p\right)\neq0$
for all $p$. It is easy to see that $g=1/H$ is a solution of \eqref{lxhg}.
So, according to Eq. \eqref{gheta} of Proposition \ref{modcon},
$X_{H}$ is the Reeb vector field of the contact form $\eta/H$, and
the next result follows form the last part of Theorem \ref{s0pi}. 
\begin{thm}
\label{impo2} Given a contact manifold $\left(M,\eta\right)$, a
fibration $\Pi:M\rightarrow N$ and a non-vanishing function $H:M\rightarrow\mathbb{R}$,
if we know a complete solution $\Sigma$ of the $\Pi$-HJE for $X_{H}$
such that 
\begin{equation}
\sigma_{\lambda}^{*}\mathsf{d}\left(\frac{\eta}{H}\right)=0,\;\;\;\forall\lambda\in\Lambda,\label{sdh}
\end{equation}
then the integral curves of $X_{H}$ can be constructed up to quadratures. 
\end{thm}
\begin{rem}
\label{equiv} If $H$ is a non-vanishing function and $\Sigma$ is
a complete solution, it is easy to show that the following statements
are equivalent: 
\end{rem}
\begin{itemize}
\item $\Sigma$ satisfies Eq. \eqref{sdh}.
\item $\Sigma$ is pseudo-isotropic w.r.t. the contact form $\eta/H$. 
\item $\Sigma$ is $1/H$-pseudo isotropic for $H$ (w.r.t. the contact
form $\eta$).
\end{itemize}

\subsubsection{Thermodynamic processes revisited}

\label{tpr} Let us go back to Section \ref{tps}. Recall we assumed
there that $a^{0}$ is constant. The complete solution given by Eqs.
\eqref{sil} and \eqref{Sil} is $g$-pseudo-isotropic with $g=e^{f}$.
To show it, first note that 
\[
X_{H}\left(f\right)=\hat{X}\left(f\right)=-a^{0}=-\xi\left(H\right),
\]
where we have used Eqs. \eqref{XH} and \eqref{lxf} and the fact
that $\xi\left(H\right)=a^{0}$. Then $X_{H}\left(e^{f}\right)+e^{f}\,\xi\left(H\right)=0$
{[}recall Eq. \eqref{lxhg}{]}. Secondly, 
\[
\left(\sigma_{\lambda}\right)_{*,\mathbf{x}}\left(\mathsf{v}\right)=\left(\mathsf{v},\left(\lambda^{k}\,\nabla g_{k}\right)\,\left(\nabla e^{-f}\right)\cdot\mathsf{v}+\left(e^{-f}\,\left(\lambda^{k}\,\mathsf{Hess}\left(g_{k}\right)\right)-\mathsf{Hess}\left(\Phi\right)\right)\cdot\mathsf{v},\nabla\left(\Phi+\lambda^{0}\,e^{-f}\right)\cdot\mathsf{v}\right),
\]
for all $\mathsf{v}\in\mathbb{R}^{n}$, and consequently 
\[
\eta\left(\left(\sigma_{\lambda}\right)_{*,\mathbf{x}}\left(\mathsf{v}\right)\right)=\left(e^{-f}\,\left(\lambda^{k}\,\nabla g_{k}\right)-\nabla\Phi\right)\cdot\mathsf{v}+\nabla\left(\Phi+\lambda^{0}\,e^{-f}\right)\cdot\mathsf{v}=e^{-f}\,\left(\lambda^{k}\,\nabla g_{k}-\lambda^{0}\,\nabla f\right)\cdot\mathsf{v},
\]
that is to say, 
\begin{equation}
\sigma_{\lambda}^{*}\eta=e^{-f}\,\mathsf{d}\left(-\lambda^{0}\,f+\lambda^{k}\,g_{k}\right).\label{W}
\end{equation}
This clearly implies that $\sigma_{\lambda}^{*}\mathsf{d}\left(e^{f}\,\eta\right)=0$,
and the $g$-pseudo-isotropy of $\Sigma$ is proved. In particular,
Eq. \eqref{W} says that the function $W:\mathbb{R}^{n}\times\mathbb{R}^{n+1}\rightarrow\mathbb{R}$
can be taken as 
\begin{equation}
W\left(\mathbf{x},\lambda\right)\coloneqq-\lambda^{0}\,f\left(\mathbf{x}\right)+\lambda^{k}\,g_{k}\left(\mathbf{x}\right).\label{W1}
\end{equation}

\begin{rem*}
As mentioned in Remark \ref{isoa0}, if $a^{0}=0$, or equivalently
$\xi\left(H\right)=0$, we have a solution in which $f=0$. Then,
the functions $\sigma_{\lambda}$ define a pseudo-isotropic complete
solution in that case. 
\end{rem*}
On the other hand, using \eqref{ck}, we have that 
\[
\begin{array}{lll}
H\circ\Sigma\left(\mathbf{x},\lambda\right) & = & H\left(\mathbf{x},e^{-f\left(\mathbf{x}\right)}\,\left(\lambda^{k}\,\nabla g_{k}\left(\mathbf{x}\right)\right)-\nabla\Phi\left(\mathbf{x}\right),\Phi\left(\mathbf{x}\right)+\lambda^{0}\,e^{-f}\right)\\
\\
 & = & e^{-f}\,\left(a^{0}\,\lambda^{0}+{\displaystyle \lambda^{k}\,c_{k}}\right).
\end{array}
\]
This means that 
\[
\left(e^{f}H\right)\circ\Sigma\left(\mathbf{x},\lambda\right)=h\left(\lambda\right)
\]
with $h:\mathbb{R}^{n+1}\rightarrow\mathbb{R}$ given by 
\[
h\left(\lambda^{0},\lambda^{1},...,\lambda^{n}\right)=a^{0}\,\lambda^{0}+\lambda^{k}\,c_{k}.
\]

Now, let us calculate the application $\varphi:\mathbb{R}^{n}\times\mathbb{R}^{n+1}\rightarrow T^{*}\mathbb{R}^{n+1}$.
If we identify $T\mathbb{R}^{n+1}$ with $\mathbb{R}^{n+1}\times\left(\mathbb{R}\times\mathbb{R}^{n}\right)$,
given $\left(\varsigma,\mathsf{w}\right)\in\mathbb{R}\times\mathbb{R}^{n}$,
it easily follows that {[}see Eq. \eqref{W1}{]} 
\[
\left\langle \mathsf{d}W\left(\mathbf{x},\lambda\right),\left(0,\left(\varsigma,\mathsf{w}\right)\right)\right\rangle =-f\left(\mathbf{x}\right)\,\varsigma+G\left(\mathbf{x}\right)\cdot\mathsf{w}
\]
and 
\[
\Sigma_{*,\left(\mathbf{x},\lambda\right)}\left(0,\left(\varsigma,\mathsf{w}\right)\right)=\left(0,e^{-f\left(\mathbf{x}\right)}\,DG\left(\mathbf{x}\right)\cdot\mathsf{w},e^{-f\left(\mathbf{x}\right)}\,\varsigma\right),
\]
where $G=\left(g_{1},...,g_{n}\right)$. Then 
\[
\eta\left(\Sigma_{*,\left(\mathbf{x},\lambda\right)}\left(0,\left(\varsigma,\mathsf{w}\right)\right)\right)=e^{-f\left(\mathbf{x}\right)}\,\varsigma
\]
and 
\[
\left\langle \varphi\left(\mathbf{x},\lambda\right),\left(\varsigma,\mathsf{w}\right)\right\rangle =\left\langle \left(\mathsf{d}W\right)\left(\mathbf{x},\lambda\right)-\Sigma^{*}\left(e^{f}\eta\right),\left(0,\left(\varsigma,\mathsf{w}\right)\right)\right\rangle =-\left(1+f\left(\mathbf{x}\right)\right)\,\varsigma+G\left(\mathbf{x}\right)\cdot\mathsf{w},
\]
what implies that 
\[
\varphi\left(\mathbf{x},\lambda\right)=\left(-1-f\left(\mathbf{x}\right),G\left(\mathbf{x}\right)\right).
\]
Here we are identifying $T_{\lambda}^{*}\mathbb{R}^{n+1}$ with $\mathbb{R}\times\mathbb{R}^{n}$.
With all that, the equations of motion reduce to 
\[
\frac{d}{dt}W_{\lambda}\left(\mathbf{x}\left(t\right)\right)=a^{0}\,\lambda^{0}+\lambda^{k}\,c_{k}
\]
and 
\[
\frac{d}{dt}\varphi_{\lambda}\left(\mathbf{x}\left(t\right)\right)=a^{0}\,\mathsf{d}\lambda^{0}+c_{k}\,\mathsf{d}\lambda^{k},
\]
being both of them equivalent to 
\[
\frac{d}{dt}f\left(\mathbf{x}\left(t\right)\right)=a^{0}
\]
and 
\begin{equation}
\frac{d}{dt}g_{i}\left(\mathbf{x}\left(t\right)\right)=-c_{i},\;\;\;i=1,...,n.\label{fig}
\end{equation}
Since the $g_{i}$'s are functionally independent functions, it is
enough to consider Eq. \eqref{fig} in order to find $\mathbf{x}\left(t\right)$.
The latter define the integral curves of each vector field $X_{H}^{\sigma_{\lambda}}\in\mathfrak{X}\left(\mathbb{R}^{n}\right)$,
and in particular of $\hat{X}=X_{H\left|\mathcal{L}\right.}$ {[}see
Eq. \eqref{X}{]}. The fact that they can be found up to quadratures
can also be explained from the following: we are dealing with a vector
field on an $n$-dimensional manifold and with $n-1$ first integrals
(see Remark \ref{cm}).

\subsubsection{A unifier point of view}

\label{un} In this paper we are actually working with a restricted
notion of the concept of contact manifold. The more general one (see
Ref. \cite{ar}) is given by a pair $\left(M,\mathcal{H}\right)$,
where $M$ is an $\left(2n+1\right)$-dimensional manifold and $\mathcal{H}$
is a distribution on $M$ of dimension $2n$: the \textbf{contact
structure}, which satisfy a certain ``non-degeneracy condition.''
Such a condition can be stated as follows: around every point of $M$
there exists an open neighborhood $U$ and a local $1$-form $\alpha$
such that $\mathsf{Ker}\,\alpha=\mathcal{H}_{\left|U\right.}$ and
$\left(\mathsf{d}\alpha\right)^{n}\wedge\alpha\neq0$. Our notion
corresponds to the case in which the contact structure can be described
by a global $1$-form. Such subclass of contact manifolds are sometime
called \textbf{co-oriented contact manifolds}. Given two of them,
say $\left(M,\eta\right)$ and $\left(M,\eta'\right)$, it is clear
that they give rise to the same contact structure $\mathcal{H}$ if
and only if there exists a non-vanishing function $g$ such that $\eta'=g\,\eta$.
This defines an equivalence relation between co-oriented contact manifolds.
On the other hand, we can say that a vector field $X\in\mathfrak{X}\left(M\right)$
is a (global) Hamiltonian contact field w.r.t. the contact structure
$\mathcal{H}$ if it is a Hamiltonian contact field w.r.t. some global
contact form $\eta$ (if exists) defining $\mathcal{H}$. As stated
in Proposition \ref{modcon}, if $X=X_{H}$ for some contact form
$\eta$, then $X=X_{gH}$ for the equivalent contact form $g\,\eta$.

So, if we think of our contact systems as triples $\left(M,\mathcal{H},X\right)$,
where $\mathcal{H}$ is a contact structure and $X$ is a Hamiltonian
contact field w.r.t. $\mathcal{H}$, we can unify the Theorems \ref{s0pi},
\ref{impo} and \ref{impo2} as follows. 
\begin{thm}
\label{gene} If we know a contact form $\eta$ defining $\mathcal{H}$,
with Reeb vector field $\xi$, such that the Hamiltonian function
$H$ of $X$ w.r.t. $\eta$ satisfies $\xi\left(H\right)=0$, and
we know a pseudo-isotropic (w.r.t. $\eta$) complete solution of the
$\Pi$-HJE for $X$, then $X$ is integrable by quadratures. 
\end{thm}

\subsection{Comments on contact non-commutative integrability}

In the context of symplectic manifolds, the notion of non-commutative
integrability was originally presented in \cite{mf}. See also \cite{j}
and references therein. For contact manifolds we shall consider the
following notion, based on Refs. \cite{jova0,jova}. Let $\left(M,\eta\right)$
be a contact manifold with Reeb vector field $\xi$. 
\begin{defn}
A contact system $\left(M,X_{H}\right)$ is \textbf{contact non-commutative
integrable} (\textbf{CNCI}) if 
\begin{enumerate}
\item $\xi\left(H\right)=0$ and 
\item there exists a fibration $F:M\rightarrow\Lambda$ such that: 
\begin{enumerate}
\item (first integrals) $\mathsf{Im}X_{H}\subset\mathsf{Ker}F_{*}$, 
\item (pseudo-isotropy) $\mathsf{Ker}F_{*}\subset\left(\mathsf{Ker}F_{*}\right)^{\bot}$, 
\item $\left(\mathsf{Ker}F_{*}\right)^{\bot}$ is integrable and 
\item $\mathsf{Im}\xi\subset\mathsf{Ker}F_{*}$, 
\end{enumerate}
\end{enumerate}
\end{defn}
If all above conditions hold, the compact and connected leaves of
$F$ are invariant tori, and action-angle-like coordinates can be
constructed for $M$ around such leaves \cite{jova0,jova}. In particular,
the system is integrable by quadratures.

\bigskip{}
 Let us explain the meaning of the items $2\left(a\right)$ and $2\left(b\right)$.
If $\Lambda$ is an open subset of $\mathbb{R}^{l}$, the item $2\left(a\right)$
says that the components $F_{1},...,F_{l}$ of $F$ must satisfy $X_{H}\left(F_{i}\right)=0$,
which means that such components are constant along the integral curves
of $X_{H}$. The item $2\left(b\right)$, on the other hand, says
that the leaves of $F$ are pseudo-isotropic submanifolds. This is
why we shall say that a submersion $F$ satisfying $2\left(a\right)$
and $2\left(b\right)$ defines a \textbf{set of pseudo-isotropic first
integrals}.

\bigskip{}

In Reference \cite{gp}, it is shown a \textit{duality} involving
(local) complete solutions of the HJE and (local) first integrals.
Concretely, it is shown that, given a complete solution $\Sigma:N\times\Lambda\rightarrow M$
of the $\Pi$-HJE for $X_{H}$, we can construct around every point
of $M$ a neighborhood $U$ and a submersion $F:U\rightarrow\Lambda$
such that 
\[
\mathsf{Im}\,X_{H\left|U\right.}\subset\mathsf{Ker}\,F_{*}\;\;\;\textrm{and}\;\;\;TU=\mathsf{Ker}\left(\Pi_{\left|U\right.}\right)_{*}\oplus\mathsf{Ker}\,F_{*}.
\]
$U$ and $F$ are given by the formulae 
\begin{equation}
U\coloneqq\Sigma\left(V\right)\;\;\;\textrm{and \;\;\;}F\coloneqq p_{\Lambda\left|U\right.}\circ\left(\Sigma_{\left|V\right.}\right)^{-1},\label{cf}
\end{equation}
where $V\subseteq N\times\Lambda$ is an open subset for which $\Sigma_{\left|V\right.}$
is a diffeomorphism with its image. Note that 
\begin{equation}
\mathsf{Ker}\,F_{*,\Sigma\left(p,\lambda\right)}=\mathsf{Im}\left(\sigma_{\lambda}\right)_{*,p},\label{kfis}
\end{equation}
for all $\left(p,\lambda\right)\in V$. This means that, if $\Sigma$
is pseudo-isotropic, the same is true for the leaves of $F$. Reciprocally,
from a submersion $F:M\rightarrow\Lambda$ satisfying 
\[
\mathsf{Im}\,X_{H}\subset\mathsf{Ker}\,F_{*}\;\;\;\textrm{and}\;\;\;TM=\mathsf{Ker}\,\Pi_{*}\oplus\mathsf{Ker}\,F_{*},
\]
we can construct, around every point of $M$, a neighborhood $U$
and a local complete solution $\Sigma$ of the $\Pi$-HJE. The involved
subset $U$ is one for which $\left(\Pi,F\right)_{\left|U\right.}$
is a diffeomorphism with its image, and $\Sigma$ is given by 
\begin{equation}
\Sigma=\left(\Pi,F\right)_{\left|U\right.}^{-1}:\Pi\left(U\right)\times F\left(U\right)\rightarrow U.\label{fc}
\end{equation}
Also, it is easy to see that \eqref{kfis} holds for all $\left(p,\lambda\right)\in\Pi\left(U\right)\times F\left(U\right)$.
Accordingly, if the leaves of $F$ are pseudo-isotropic, so is $\Sigma$.

Summarizing, a pseudo-isotropic complete solution gives rise to a
set of local pseudo-isotropic first integrals \textit{via} Eq. \eqref{cf},
and a set of pseudo-isotropic first integrals gives rise to a local
pseudo-isotropic complete solution \textit{via} Eq. \eqref{fc}.

\bigskip{}

From this duality, given a contact system $\left(M,X_{H}\right)$
such that $\xi\left(H\right)=0$, if we know a set of pseudo-isotropic
first integrals $F$, i.e. a submersion satisfying the items $2\left(a\right)$
and $2\left(b\right)$ of the last definition, then, according to
Theorem \ref{s0pi}, such a system is integrable up to quadratures.
In particular, the other two items {[}$2\left(c\right)$ and $2\left(d\right)${]}
are not needed to ensure integrability by quadratures! 
\begin{rem*}
An analogous result was proved in Ref. \cite{gp}, but in the context
of symplectic and Poisson manifolds (see also Ref. \cite{g}).
\end{rem*}
Moreover, if we consider the notions of contact manifolds and contact
systems given in Section \ref{un}, from Theorem \ref{gene} and the
above discussion we have the next result. 
\begin{thm}
Given a contact system $\left(M,\mathcal{H},X\right)$, if we know
a contact form $\eta$ defining $\mathcal{H}$, with Reeb vector field
$\xi$, such that the Hamiltonian function $H$ of $X$ w.r.t. $\eta$
satisfies $\xi\left(H\right)=0$, and we know a set of pseudo-isotropic
first integrals $F$ (for the contact form $\eta$), then $X$ is
integrable by quadratures. 
\end{thm}
\bigskip{}

As an example, consider the system of Section \ref{tp}. It can be
described as a triple $\left(M,\mathcal{H},X\right)$ with $M=\mathbb{R}^{2n+1}$,
$\mathcal{H}=\mathsf{Ker}\,\eta$, $\eta\coloneqq y_{i}\,\mathsf{d}x^{i}+\mathsf{d}z$
and $X$ given by \eqref{XH}. We have shown in Section \ref{tpr}
that the complete solution $\Sigma$ given by \eqref{sil} and \eqref{Sil}
is $e^{f}$-pseudo-isotropic. So, defining $\eta'=e^{f}\eta$ and
calling $\xi'$ its Reeb vector field, form the Proposition \ref{modcon}
(and the comment below it) we have that the Hamiltonian function of
$X$ w.r.t. $\eta'$ is $H'=e^{f}H$ {[}see Eq. \eqref{Htp}{]} and
satisfies $\xi'\left(H'\right)=0$. Also, we have that $\Sigma$ is
pseudo-isotropic w.r.t. $\eta'$ and that $\mathcal{H}=\mathsf{Ker}\,\eta'$.
Now, let us find the first integrals related to $\Sigma$. It is easy
to see that 
\[
\Sigma^{-1}\left(\mathbf{x},\mathbf{y},z\right)=\left(\mathbf{x},e^{f\left(\mathbf{x}\right)}\,\left(z-\Phi\left(\mathbf{x}\right)\right),e^{f\left(\mathbf{x}\right)}\,\left(\mathbf{y}-\nabla\Phi\left(\mathbf{x}\right)\right)^{t}\cdot\left[DG\left(\mathbf{x}\right)\right]^{-1}\right),
\]
where $DG\left(\mathbf{x}\right)$ is the differential of the function
$G\coloneqq\left(g_{1},...,g_{n}\right):\mathbb{R}^{n}\rightarrow\mathbb{R}^{n}$.
Then, Eq. \eqref{cf} (and the duality above discussed) ensures that
\[
F\left(\mathbf{x},\mathbf{y},z\right)\coloneqq\left(F_{0},...,F_{n}\right)\left(\mathbf{x},\mathbf{y},z\right)\coloneqq\left(e^{f\left(\mathbf{x}\right)}\,\left(z-\Phi\left(\mathbf{x}\right)\right),e^{f\left(\mathbf{x}\right)}\,\left(\mathbf{y}-\nabla\Phi\left(\mathbf{x}\right)\right)^{t}\cdot\left[DG\left(\mathbf{x}\right)\right]^{-1}\right)
\]
defines a set of $n+1$ pseudo-isotropic first integrals of $X$ (w.r.t.
$\eta'$). So we have a contact form and a fibration satisfying the
theorem above. 
\begin{rem*}
Using Eq. \eqref{sphp}, it is easy to show that $\xi'\left(F_{0}\right)=1$.
This means that $\mathsf{Im}\,\xi'\nsubseteq\mathsf{Ker}\,F_{*}$.
Thus, the item $2\left(d\right)$ of above definition does not hold
for $F$. Also, it can be shown that $\left(\mathsf{Ker}\,F_{*}\right)^{\bot}$,
which in this case is equal to $\mathsf{Ker}\,F_{*}+\left\langle \xi'\right\rangle $,
is not integrable, i.e. item $2\left(c\right)$ does not hold either. 
\end{rem*}

\section{Another integration procedure}

As we have seen in Section \ref{integ}, if we have a contact system
$\left(M,X_{H}\right)$ such that $\xi\left(H\right)=0$, and we know
a pseudo-isotropic complete solution for it, then $X_{H}$ can be
integrated up to quadratures. If condition $\xi\left(H\right)=0$
does not hold, to ensure integrability up to quadratures we can try
to find a solution $g$ of Eq. \eqref{lxhg} and a $g$-pseudo-isotropic
complete solution. When $H$ is a non-vanishing function, we can take
$g=1/H$. Of course, we could also do that in the open submanifold
where $H\neq0$. This is good because, in such a case, we do not need
to solve Eq. \eqref{lxhg}. But, what can we do on the subset of $M$
where $\xi\left(H\right)\neq0$ and $H=0$? Do we have an alternative
to the $g$-pseudo-isotropic complete solutions? To give an answer
to these questions, we first need some definitions.

\subsection{The $\xi\left(H\right)\protect\neq0$ case and the bi-isotropy condition}

Let $\left(M,X_{H}\right)$ be a contact system such that $\xi\left(H\right)\neq0$.
For simplicity, suppose that 
\[
M_{0}\coloneqq H^{-1}\left(0\right)=\left\{ m\in M:H\left(m\right)=0\right\} 
\]
is a closed regular submanifold of $M$ (i.e. $0$ is a regular value
of $H$) and denote by $M_{1}$ the (open) complement of $M_{0}$,
i.e. 
\begin{equation}
M_{1}\coloneqq\left(H^{-1}\left(0\right)\right)^{c}=\left\{ m\in M':H\left(m\right)\neq0\right\} .\label{m1}
\end{equation}
Consider the related inclusions $\mathfrak{i}_{0,1}:M_{0,1}\rightarrow M$. 
\begin{prop}
\label{12} Under above conditions and notation, 
\begin{enumerate}
\item the pair $\left(M_{0},\mathsf{d}\eta_{0}\right)$, with $\eta_{0}\coloneqq\mathfrak{i}_{0}^{*}\eta$,
is a symplectic manifold. Also, $M_{0}$ is an $X_{H}$-invariant
submanifold. In particular, there exists a vector field $Y\in\mathfrak{X}\left(M_{0}\right)$
such that 
\begin{equation}
\mathfrak{i}_{0*}\circ Y=X_{H}\circ\mathfrak{i}_{0};\label{irel}
\end{equation}
\item the pair $\left(M_{1},\eta_{1}\right)$, with $\eta_{1}\coloneqq\mathfrak{i}_{1}^{*}\eta$,
is a contact manifold and, defining $H_{1}\coloneqq\mathfrak{i}_{1}^{*}H$,
\[
\mathfrak{i}_{1*}\circ X_{H_{1}}=X_{H}\circ\mathfrak{i}_{1}.
\]
\end{enumerate}
\end{prop}
\textit{Proof}. $1.$ The condition $\xi\left(H\right)\neq0$ is equivalent
to 
\begin{equation}
\xi\left(p\right)\notin\mathsf{Ker}\left(\mathsf{d}H\left(p\right)\right)=T_{p}M_{0},\;\;\;\forall p\in M_{0}.\label{kdh}
\end{equation}
So, since $\mathsf{Ker}\left(\mathsf{d}\eta\right)_{p}=\left\langle \xi\left(p\right)\right\rangle $,
it follows that $\mathfrak{i}_{0}^{*}\mathsf{d}\eta$ is a symplectic
form on $M_{0}$. On the other hand, since $X_{H}\left(H\right)=H\,\xi\left(H\right)$
{[}see the first part of Eq. \eqref{dcvf}{]}, we have that 
\[
X_{H}\left(p\right)\in\mathsf{Ker}\left(\mathsf{d}H\left(p\right)\right)=T_{p}M_{0},\;\;\;\forall p\in M_{0},
\]
what says precisely that $M_{0}$ is $X_{H}$-invariant.

\bigskip{}

$2.$ Since $M_{1}$ is a open submanifold, it is clear that $\mathfrak{i}_{1}^{*}\eta$
is a contact form. For the same reason, if we restrict \eqref{drvf}
to a point $p\in M_{1}$, we shall find that the Reeb vector field
$\xi_{1}$ of $\left(M_{1},\eta_{1}\right)$ is given $\xi_{1}\left(p\right)=\xi\left(p\right)$
(i.e. the restriction of $\xi$ to $M_{1}$). Moreover, if we make
the same restriction on Eq. \eqref{dcvf}, we shall have that 
\[
\mathsf{d}\eta\left(X_{H}\left(p\right),v\right)=\left\langle -\mathsf{d}H_{1}\left(p\right)+\xi_{1}\left(H_{1}\right)\left(p\right)\,\left(\eta_{1}\right)_{p},v\right\rangle ,\;\;\;\forall v\in T_{p}M_{1},
\]
and 
\[
\eta\left(X_{H}\left(p\right)\right)=H_{1}\left(p\right).
\]
This means, omitting the inclusion $\mathfrak{i}_{1}$, that $X_{H}\left(p\right)=X_{H_{1}}\left(p\right)$
for all $p\in M_{1}$, as we wanted to show. $\;\;\;\triangle$

\bigskip{}

It is clear that, if we have the integral curves of the vector fields
$Y=\left(X_{H}\right)_{\left|M_{0}\right.}\in\mathfrak{X}\left(M_{0}\right)$
and $X_{H_{1}}=\left(X_{H}\right)_{\left|M_{1}\right.}\in\mathfrak{X}\left(M_{1}\right)$
(see proposition above), then we have those of $X_{H}$. So, if we
want to find such curves by means of a complete solution, we need:
\begin{enumerate}
\item to restrict, in some sense, the complete solutions to the $X$-invariant
submanifolds $M_{0}$ and $M_{1}$;
\item to impose on the resulting restrictions appropriate conditions to
ensure integrability by quadratures.
\end{enumerate}
\bigskip{}

For the first point, we have the next definition (based on Definition
4.17 of Ref. \cite{gp}). Consider an arbitrary manifold $M$, a vector
field $X\in\mathfrak{X}\left(M\right)$, a fibration $\Pi:M\rightarrow N$
and a complete solution $\Sigma:N\times\Lambda\rightarrow M$.
\begin{defn}
\label{grc}Given an $X$-invariant submanifold $M'\subset M$, we
shall say that\textbf{ }$\Sigma$\textbf{ restricts to }$M'$ if there
exist submanifolds $N',\Lambda'$ of $N$ and $\Lambda$, respectively
such that $\Pi_{|M'}:M'\rightarrow N'$ is a fibration, $\Sigma\left(N'\times{\Lambda}'\right)\subset M'$
and 
\[
\Sigma_{M'}:=\Sigma_{|N'\times\Lambda'}:N'\times\Lambda'\rightarrow M'
\]
is a complete solution of the $\Pi_{\left|M'\right.}$-HJE for $X_{\left|M'\right.}$.
We shall call $\Sigma_{M'}$ a \textbf{restriction} of $\Sigma$ to
$M'$. On the other hand, if for every $p\in M'$ there exists an
open neighborhood $U\subseteq M'$ of $p$ such that $\Sigma$\textbf{
restricts to }$U$, then we shall say that $\Sigma$\textbf{ locally
restricts to }$M'$, and we shall call \textbf{local restriction}
of $\Sigma$ to $M'$ any related map $\Sigma_{U}$.
\end{defn}
\begin{rem}
If $M'$ is a open submanifold, it is always true that $M$ is $X$-invariant
and that $\Sigma$ locally restricts to $M'.$
\end{rem}
For the second point, we shall consider the following conditions.
Let us go back to the contact context.
\begin{defn}
\label{bis}Under the conditions of Proposition \ref{12}, we shall
say that a complete solution $\Sigma$ is \textbf{bi-isotropic} if: 
\begin{enumerate}
\item $\Sigma$ restricts to $M_{0}$ and for some restriction $\Sigma_{0}\coloneqq\Sigma_{M_{0}}$
we have that
\begin{equation}
\left(\sigma_{0\lambda}\right)^{*}\eta_{0}=0\label{ic}
\end{equation}
for all its partial solutions $\sigma_{0\lambda}$.
\item $\Sigma$ restricts to $M_{1}$ and for some restriction $\Sigma_{1}\coloneqq\Sigma_{M_{1}}$
we have that 
\begin{equation}
\left(\sigma_{1\lambda}\right)^{*}\mathsf{d}\left(\frac{\eta_{1}}{H_{1}}\right)=0\label{pic}
\end{equation}
for all its partial solutions $\sigma_{1\lambda}$.
\end{enumerate}
The \textbf{locally bi}-\textbf{isotropic} version consists in replacing
above restrictions $\Sigma_{0}$ and $\Sigma_{1}$ by local restrictions
around every point of the corresponding manifolds.
\end{defn}
\begin{rem}
\label{ssm} Note that Eq. \eqref{ic} implies that $\left(\sigma_{0\lambda}\right)^{*}\mathsf{d}\eta_{0}=0$,
i.e. $\Sigma_{0}$ is isotropic in the sense of symplectic manifolds
(see Ref. \cite{gp}). On the other hand, Eq. \eqref{ic} says exactly
that $\Sigma_{1}$ is $1/H_{1}$-pseudo-isotropic (w.r.t. the contact
form $\eta_{1}$) or, equivalently, pseudo-isotropic w.r.t. $\eta_{1}/H_{1}$
(see Remark \ref{equiv}).
\end{rem}

\subsection{The integration process for the $\xi\left(H\right)\protect\neq0$
case }

\label{ishnn} Let $\Sigma$ be a bi-isotropic complete solution.
Then, in particular, we have a complete solution $\Sigma_{1}:N_{1}\times\hat{\Lambda}\rightarrow M_{1}$
of the $\Pi_{\left|M_{1}\right.}$-HJE for the vector field $\left(X_{H}\right)_{\left|M_{1}\right.}=X_{H_{1}}$,
which is pseudo-isotropic w.r.t. $\eta_{1}/H_{1}$. According to Theorem
\ref{impo2}, this implies that $X_{H_{1}}$ is integrable up to quadratures.
So, the vector field $X_{H}$ can be explicitly integrated (up to
quadratures) along $M_{1}$. Let us see what happens along $M_{0}$.

Bi-isotropy condition also ensures that we have a complete solution
$\Sigma_{0}:N_{0}\times\hat{\Lambda}\rightarrow M_{0}$ of the $\Pi_{\left|M_{0}\right.}$-HJE
for the vector field $\left(X_{H}\right)_{\left|M_{0}\right.}=Y\in\mathfrak{X}\left(M_{0}\right)$.
Thus, $\Sigma_{0}$ satisfies the identity 
\begin{equation}
\Sigma_{0*}\circ Y^{\Sigma_{0}}=Y\circ\Sigma_{0}.\label{sys}
\end{equation}
Also, Eq. \eqref{ic} ensures that
\begin{equation}
i_{Y^{\Sigma_{0}}}\Sigma_{0}^{*}\eta_{0}=0.\label{yso}
\end{equation}
{[}This identity could also be deduced from the second part of \eqref{ccont}{]}.
In addition, using the first part of \eqref{ccont} and the Eqs. \eqref{irel}
and \eqref{sys}, for every vector $v\in T_{p}N_{0}\times T_{\lambda}\hat{\Lambda}$
we have that 
\[
\begin{array}{lll}
\left(i_{Y^{\Sigma_{0}}}\Sigma_{0}^{*}\mathsf{d}\eta_{0}\right)_{\left(p,\lambda\right)}\left(v\right) & = & \left(\mathsf{d}\eta\right)_{\Sigma_{0}\left(p,\lambda\right)}\left(\mathfrak{i}_{0*}\circ\Sigma_{0*}\circ Y^{\Sigma_{0}}\left(p,\lambda\right),\left(\mathfrak{i}_{0}\circ\Sigma_{0}\right)_{*,\left(p,\lambda\right)}\left(v\right)\right)\\
 & = & \left(\mathsf{d}\eta\right)_{\Sigma_{0}\left(p,\lambda\right)}\left(X_{H}\circ\mathfrak{i}_{0}\circ\Sigma_{0}\left(p,\lambda\right),\left(\mathfrak{i}_{0}\circ\Sigma_{0}\right)_{*,\left(p,\lambda\right)}\left(v\right)\right)\\
 & = & \left\langle \left(\xi\left(H\right)\,\eta-dH\right)\circ\mathfrak{i}_{0}\circ\Sigma_{0}\left(p,\lambda\right),\left(\mathfrak{i}_{0}\circ\Sigma_{0}\right)_{*,\left(p,\lambda\right)}\left(v\right)\right\rangle ,
\end{array}
\]
and consequently, 
\[
i_{Y^{\Sigma_{0}}}\Sigma_{0}^{*}\mathsf{d}\eta_{0}=\Sigma_{0}^{*}\mathfrak{i}_{0}^{*}\left(\xi\left(H\right)\,\eta-dH\right)=\Sigma_{0}^{*}\left(\left(\xi\left(H\right)\circ\mathfrak{i}_{0}\right)\,\eta_{0}\right),
\]
where we have used that $\mathfrak{i}_{0}^{*}\mathsf{d}H=0$. Combining
\eqref{yso} and the last equation, we have that 
\begin{equation}
\mathsf{L}_{Y^{\Sigma_{0}}}\Sigma_{0}^{*}\eta_{0}=\varsigma\,\Sigma_{0}^{*}\eta_{0},\label{prop}
\end{equation}
with
\begin{equation}
\varsigma\coloneqq\xi\left(H\right)\circ\mathfrak{i}_{0}\circ\Sigma_{0}.\label{seda}
\end{equation}
From Eq. \eqref{prop} we shall construct the integral curves of $Y$.
To do that, as in previous sections, consider the function 
\[
\hat{\varphi}:N_{0}\times\hat{\Lambda}\rightarrow T^{*}\hat{\Lambda}
\]
such that $\hat{\varphi}\left(p,\lambda\right)\in T_{\lambda}^{*}\hat{\Lambda}$,
for all $\left(p,\lambda\right)\in N_{0}\times\hat{\Lambda}$, and
\begin{equation}
\left\langle \hat{\varphi}\left(p,\lambda\right),z\right\rangle =-\left\langle \left(\Sigma_{0}^{*}\eta_{0}\right)\left(p,\lambda\right),\left(0,z\right)\right\rangle ,\;\;\;\forall z\in T_{\lambda}\hat{\Lambda},\label{defi-1}
\end{equation}
with its related functions $\hat{\varphi}_{\lambda}:N_{0}\rightarrow T_{\lambda}^{*}\hat{\Lambda}$.
Using the isotropy property of $\Sigma_{0}$ in the sense of symplectic
manifolds (see Remark \ref{ssm}), the next result can be shown as
in Ref. \cite{gp} (see Proposition 3.16 there). 
\begin{prop}
\label{imm3} Each function $\hat{\varphi}_{\lambda}$ (and consequently
$\hat{\varphi}$) is an immersion. 
\end{prop}
Let us go back to Eq. \eqref{prop}. Given an integral curve $\left(\gamma\left(t\right),\lambda\right)$
of $Y^{\Sigma_{0}}$, in terms of $\hat{\varphi}_{\lambda}$ such
an equation reads 
\begin{equation}
\frac{d}{dt}\hat{\varphi}_{\lambda}\left(\gamma\left(t\right)\right)=\varsigma_{\lambda}\left(\gamma\left(t\right)\right)\,\hat{\varphi}_{\lambda}\left(\gamma\left(t\right)\right),\label{pceq}
\end{equation}
with $\varsigma_{\lambda}\left(\gamma\left(t\right)\right)\coloneqq\varsigma\left(\gamma\left(t\right),\lambda\right)$.
{[}This can be proved by using Eq. \eqref{conm}, but replacing $X_{H}^{\Sigma}$,
$\eta$ and $W$ by $Y^{\Sigma_{0}}$, $\eta_{0}$ and $0$, respectively,
and using property \eqref{ic}{]}. On the other hand, last proposition
implies that, around every point of $N_{0}$, and fixing a basis of
$T_{\lambda}^{*}\hat{\Lambda}$, there exists a number $k=\mathsf{dim}N_{0}$
of coordinates of $\hat{\varphi}_{\lambda}$ (w.r.t. that basis) that
defines a coordinate system $\psi=\left(\psi_{1},...,\psi_{k}\right)$
for $N_{0}$. Then, in terms of $\psi$, equation above translates
to 
\begin{equation}
\dot{\psi}_{i}\left(t\right)=\varsigma_{\lambda}\left(\psi^{-1}\left(\psi_{1}\left(t\right),...,\psi_{k}\left(t\right)\right)\right)\,\psi_{i}\left(t\right),\label{ceq}
\end{equation}
where for simplicity we are writing $\psi_{i}\left(\gamma\left(t\right)\right)=\psi_{i}\left(t\right)$.
Let us show that above equations can be solved by quadratures. Note
first that, since $\xi\left(H\right)\neq0$, the function $\varsigma_{\lambda}$
is non-vanishing. Then, beside the trivial solution, we can look for
a solution where, for instance, $\psi_{1}\left(t\right)$ is not identically
zero. For such a solution we would have that 
\[
\dot{\psi}_{i}\left(t\right)=\frac{\dot{\psi}_{1}\left(t\right)}{\psi_{1}\left(t\right)}\,\psi_{i}\left(t\right),
\]
and consequently $\psi_{i}\left(t\right)=c_{i}\,\psi_{1}\left(t\right)$
for some constant $c_{i}$. Replacing this in \eqref{ceq} for $i=1$,
we obtain the equation 
\[
\dot{\psi}_{1}\left(t\right)=F\left(\psi_{1}\left(t\right)\right)\;\;\;\textrm{with}\;\;\;F\left(x\right)=\varsigma_{\lambda}\left(\psi^{-1}\left(x,c_{2}\,x,...,c_{k}\,x\right)\right)\,x.
\]
It is clear that if we find a function $\psi_{1}$ satisfying the
equation above, which can be constructed up to one quadrature, we
have the solution we are looking for. Concluding, we have the next
result. 
\begin{thm}
\label{s1pi}Given a contact manifold $\left(M,\eta\right)$, a fibration
$\Pi:M\rightarrow N$ and a function $H:M\rightarrow\mathbb{R}$ such
that $\xi\left(H\right)$ is a non-vanishing function, if we know
a bi-isotropic complete solution $\Sigma$ of the $\Pi$-HJE for $X_{H}$,
then the integral curves of $X_{H}$ can be constructed up to quadratures. 
\end{thm}
\begin{rem*}
Of course, if instead of a bi-isotropic complete solution on $M$
we have a complete solution satisfying Eq. \eqref{ic} on $M_{0}$
and a pseudo-isotropic complete solution on $M_{1}$ (independent
of each other), then we can proceed exactly as we did above. We have
introduced the concept of bi-isotropy just with the aim of presenting
a set of sufficient conditions (that ensure integrability by quadratures)
involving the whole phase space $M$ of the system.
\end{rem*}

\subsection{An illustrative example}

Consider the manifold $\mathbb{R}^{3}$, denote its global coordinates
as $\left(q,p,s\right)$, and define on it the function $H$ such
that 
\[
H\left(q,p,s\right)=\frac{p^{2}+q^{2}}{2}-\alpha s,
\]
with $\alpha\neq0$. We shall consider on $\mathbb{R}^{3}$ the contact
form $\eta=p\,\mathsf{d}q+\mathsf{d}s$. The related contact system,
with 
\[
X_{H}\left(q,p,s\right)=p\,\frac{\partial}{\partial q}-\left(q+\alpha\,p\right)\,\frac{\partial}{\partial q}+\left(\frac{q^{2}-p^{2}}{2}-\alpha\,s\right)\,\frac{\partial}{\partial s},
\]
is sometimes known as the $1$-\textit{dimensional damped oscillator}. 
\begin{rem*}
This system has also been studied in Ref. \cite{ds} (see also Ref.
\cite{bra}). Here, besides finding a complete solution of its $\Pi$-HJE,
for a given fibration $\Pi$ (see below), we also find an expression
for some of its trajectories. 
\end{rem*}
Note that 
\begin{equation}
\xi\left(H\right)=\frac{\partial}{\partial s}H=-\alpha\neq0.\label{sedaH}
\end{equation}
Consider also the fibration $\Pi:\mathbb{R}^{3}\rightarrow\mathbb{R}$
such that $\Pi\left(q,p,s\right)=q$. A section 
\[
\sigma:\mathbb{R}\rightarrow\mathbb{R}^{3}\;:\;q\mapsto\left(q,\phi\left(q\right),\chi\left(q\right)\right)
\]
is a solution of the $\Pi$-HJE for above system if and only if 
\begin{equation}
\phi'\left(q\right)\,\phi\left(q\right)=-q-\alpha\,\phi\left(q\right)\label{sig1}
\end{equation}
and 
\begin{equation}
\chi'\left(q\right)\,\phi\left(q\right)=\frac{q^{2}-\left(\phi\left(q\right)\right)^{2}}{2}-\alpha\,\chi\left(q\right).\label{sig2}
\end{equation}
Since the first equation is homogeneous, a solution $\phi$ can be
found up to quadratures, unless locally. In fact, for $\alpha\neq2$,
the general solution of \eqref{sig1} is given by the algebraic equation
\begin{equation}
\frac{\left(\phi\left(q\right)-a_{+}\,q\right)^{b_{+}}}{\left(\phi\left(q\right)-a_{-}\,q\right)^{b_{-}}}=\frac{\lambda_{1}}{\sqrt{q}},\label{gs}
\end{equation}
with $a_{\pm}=\left(-\alpha/2\right)\pm\sqrt{1-\left(\alpha/2\right)^{2}}$,
$b_{\pm}=a_{\pm}/\sqrt{4-\alpha^{2}}$ and $\lambda_{1}$ arbitrary.
For $\alpha=2$, the general solution is given by 
\[
\left(\phi\left(q\right)+q\right)\,\exp\left(\frac{q}{\phi\left(q\right)+q}\right)=\lambda_{1}.
\]
In any case, given a solution $\phi$ of \eqref{sig1}, the general
solution of \eqref{sig2} is 
\begin{equation}
\chi\left(q\right)=\frac{q^{2}+\left(\phi\left(q\right)\right)^{2}}{2\alpha}+\lambda_{2}\,\exp\left(\int\frac{-1}{\phi\left(r\right)}dr\right),\label{sol2}
\end{equation}
with $\lambda_{2}$ also arbitrary. Taking into account the dependence
on $\lambda_{1}$ and $\lambda_{2}$ of above solutions, we shall
refer to them as $\phi^{\lambda_{1}}$ and $\chi^{\lambda_{1},\lambda_{2}}$.
It can be shown that the function $\Sigma:N\times\Lambda\rightarrow\mathbb{R}^{3}$,
with $N=\mathbb{R}-\left\{ 0\right\} $ and $\Lambda=\left(\mathbb{R}-\left\{ 0\right\} \right)\times\mathbb{R}$,
such that 
\[
\Sigma\left(q,\lambda_{1},\lambda_{2}\right)\coloneqq\left(q,\phi^{\lambda_{1}}\left(q\right),\chi^{\lambda_{1},\lambda_{2}}\left(q\right)\right),
\]
is a global diffeomorphism with its (open dense) image $M\coloneqq\Sigma\left(N\times\Lambda\right)\subseteq\mathbb{R}^{3}$,
and consequently a complete solution along $M$. Moreover, since for
$\lambda_{2}=0$ we have that {[}see \eqref{sol2}{]} 
\[
H\left(q,\phi^{\lambda_{1}}\left(q\right),\chi^{\lambda_{1},0}\left(q\right)\right)=\frac{q^{2}+\left(\phi^{\lambda_{1}}\left(q\right)\right)^{2}}{2}-\alpha\,\chi^{\lambda_{1},0}\left(q\right)=0,
\]
then $\Sigma$ restricts to $M_{0}\coloneqq H^{-1}\left(0\right)\cap M$
(which is a closed regular submanifold of $M$). Its restriction $\Sigma_{0}$
is defined on $N\times\hat{\Lambda}$, with $\hat{\Lambda}=\left(\mathbb{R}-\left\{ 0\right\} \right)\times\left\{ 0\right\} \simeq\mathbb{R}-\left\{ 0\right\} $.
Also, related to $M_{\text{1}}\coloneqq\left(H^{-1}\left(0\right)\right)^{c}\cap M$,
we have an open restriction $\Sigma_{1}$ with $\Lambda'=\left(\mathbb{R}-\left\{ 0\right\} \right)\times\left(\mathbb{R}-\left\{ 0\right\} \right)$.
(Recall Definition \ref{grc}). For dimensional reasons, both $\Sigma_{0}$
and $\Sigma_{1}$ satisfy the conditions of Definition \ref{bis}.
Let us concentrate on $M_{0}$. According to Eq. \eqref{defi-1},
\[
\hat{\varphi}_{\lambda_{1}}:\mathbb{R}-\left\{ 0\right\} \rightarrow T_{\lambda_{1}}^{*}\left(\mathbb{R}-\left\{ 0\right\} \right)
\]
is given by 
\[
\hat{\varphi}_{\lambda_{1}}\left(q\right)=\frac{1}{2}\,\frac{\partial}{\partial\lambda_{1}}\left(\phi^{\lambda_{1}}\left(q\right)\right)^{2}\,\mathsf{d}\lambda_{1}.
\]
Then, integrating \eqref{pceq} for a curve $q\left(t\right)$ on
$N$ (taking into account that in this case $\varsigma=-\alpha$ -see
\eqref{seda} and \eqref{sedaH}-), we have that 
\[
\frac{\partial}{\partial\lambda_{1}}\left(\phi^{\lambda_{1}}\left(q\left(t\right)\right)\right)^{2}=\frac{\partial}{\partial\lambda_{1}}\left(\phi^{\lambda_{1}}\left(q\left(0\right)\right)\right)^{2}\,e^{-\alpha t}.
\]
So, in order to find the integral curves of $X_{H}$ along $M_{0}$,
it only rests to solve the above algebraic equation for $q\left(t\right)$.
Once this equation is solved, such integral curves can be written
as 
\[
\left(q\left(t\right),\phi^{\lambda_{1}}\left(q\left(t\right)\right),\frac{\left(q\left(t\right)\right)^{2}+\left(\phi^{\lambda_{1}}\left(q\left(t\right)\right)\right)^{2}}{2\alpha}\right).
\]

\subsection{Addressing the general situation}

Consider a contact Hamiltonian system $\left(M,X_{H}\right)$. In
order to ensure the integrability up to quadratures of its equations
of motion, we can look for a $g$-pseudo-isotropic solution, as we
have seen in Section \ref{ges}, or we can combine the results of
Sections \ref{si0} and \ref{ishnn}. For the second option, define
\[
M'\coloneqq\left\{ m\in M:\xi\left(H\right)\neq0\right\} \;\;\;\textrm{and}\;\;\;M''\coloneqq M\left\backslash M'\right.=\left\{ m\in M:\xi\left(H\right)=0\right\} .
\]
Inside $M'$, consider the subsets 
\[
M_{0}\coloneqq H^{-1}\left(0\right)\cap M'\;\;\;\textrm{and}\;\;\;M_{1}\coloneqq M'\left\backslash M_{0}\right.,
\]
and inside $M''$ consider its interior $M_{2}\coloneqq\mathsf{int}\,M''$.
Note that on the open submanifold $M'$ we have, as in the previous
sections, a symplectic manifold $\left(M_{0},\eta_{0}\right)$ and
a contact manifold $\left(M_{1},\eta_{1}\right)$. In the latter,
we have that $\left(X_{H}\right)_{\left|M_{1}\right.}=X_{H_{\left|M_{1}\right.}}$.
Also, since $M_{2}$ is an open submanifold, we have another contact
manifold $\left(M_{2},\eta_{2}\right)$ where $\left(X_{H}\right)_{\left|M_{2}\right.}=X_{H_{\left|M_{2}\right.}}$.
With these observations in mind, we can enunciate the following theorem.
Its proof is contained in those of Theorems \ref{s0pi} and \ref{s1pi}. 
\begin{thm}
Given a contact manifold $\left(M,\eta\right)$, a fibration $\Pi:M\rightarrow N$
and a function $H:M\rightarrow\mathbb{R}$, if we know: 
\begin{enumerate}
\item a bi-isotropic solution on $M_{0}\cup M_{1}=M'$; 
\item a pseudo-isotropic solution on $M_{2}=\mathsf{int}\,M''$. 
\end{enumerate}
then the integral curves of $X_{H}$ can be constructed up to quadratures
on the open dense subset $M\left\backslash \partial M''\right.$. 
\end{thm}
Since the topological boundary $\partial M''$ of $M''$ is not in
general an invariant submanifold, the trajectories with initial conditions
inside such a subset must be studied separately (case by case).

\section*{Acknowledgements}

This work was partially supported by Ministerio de Ciencia, Innovación
y Universidades (Spain), grants MTM2015-64166-C2-2-P and PGC2018-098265-B-C32
(E.P.), the European Community IRSES-project GEOMECH-246981 (S.G.
and E.P.) and CONICET (S.G.).

\end{document}